\documentclass[11pt]{amsart}
\usepackage{graphicx}
\usepackage{amssymb}
\usepackage{xcolor}
\date{}

\newtheorem{theorem}{Theorem}
\newtheorem*{remark}{Remark}
\newtheorem{corollary}{Corollary}
\newtheorem{definition}{Definition}
\newtheorem*{ex}{Example}
\newtheorem{lemma}{Lemma}

\renewcommand{\S}{\mathbb{S}}

\newcommand{\D}{\mathbb{D}}

\newcommand{\Z}{\mathbb{Z}}

\newcommand{\R}{\mathbb{R}}
\newcommand{\C}{\mathbb{C}}

\newcommand{\SC}{\mathcal{C}}

\newcommand{\La}{\Lambda}
\newcommand{\la}{\lambda}

\newcommand{\dd}{\partial}

\newcommand{\sse}{\subset}
\newcommand{\lr}{\rightarrow}

\newcommand{\st}{\text{st}}

\begin{document}
	\title{Weave Realizability for $D-$type}
	\author{James Hughes}
	\maketitle


	

	
	
	
	
	

	
	\begin{center}
		\begin{abstract}
			We study exact Lagrangian fillings of Legendrian links of $D_n$-type in the standard contact 3-sphere. The main result is the existence of a Lagrangian filling, represented by a weave, such that any algebraic quiver mutation of the associated intersection quiver can be realized as a geometric weave mutation. The method of proof is via Legendrian weave calculus and a construction of appropriate 1-cycles whose geometric intersections realize the required algebraic intersection numbers. In particular, we show that in $D$-type, each cluster chart of the moduli of microlocal rank-1 sheaves is induced by at least one embedded exact Lagrangian filling. Hence, the Legendrian links of $D_n$-type have at least as many Hamiltonian isotopy classes of Lagrangian fillings as cluster seeds in the $D_n$-type cluster algebra, and their geometric exchange graph for Lagrangian disk surgeries contains the cluster exchange graph of $D_n$-type.
		\end{abstract}    
	\end{center}

	\section{Introduction}
	
	Legendrian links in contact 3-manifolds \cite{Bennequin83,ArnoldSing} are central to the study of 3-dimensional contact topology \cite{OzbagciStipsicz04,Geiges08}. Recent developments \cite{CasalsZaslow, CasalsGao, CasalsNg} have revealed new phenomena regarding their Lagrangian fillings, including the existence of many
	Legendrian links  $\Lambda\subseteq(\mathbb{S}^3,\xi_{\st})$ with infinitely many (smoothly isotopic) Lagrangian fillings in the Darboux 4-ball $(\D^4,\omega_{\st})$ 
	which are not Hamiltonian isotopic. The relationship between cluster algebras and Lagrangian fillings \cite{CasalsZaslow,GSW} has also led to new conjectures on the classification of Lagrangian fillings \cite{CasalsLagSkel}. In particular, \cite[Conjecture 5.1]{CasalsLagSkel} introduced a conjectural ADE classification of Lagrangian fillings. The object of this manuscript is to study $D$-type and prove part of the conjectured classification.
	
	The $A$-type was studied in \cite{EHK,Pan-fillings}, via Floer-theoretic methods, and in \cite{STWZ,TreumannZaslow} via microlocal sheaves. Their main result is that the $A_{n}$-Legendrian link $\lambda(A_{n})\subseteq(\mathbb{S}^3,\xi_{\st})$, which is the max-tb representative of the $(2, n+1)$-torus link, has at least a Catalan number $C_{n+1}=\frac{1}{n+2}{2n+2\choose n+1}$ of embedded exact Lagrangian fillings, where $C_{n+1}$ is precisely the number of cluster seeds in the finite type $A_n$ cluster algebra \cite{FWZ2}. We will show that the same holds in $D$-type, namely that $D_n$-type Legendrian links have at least as many distinct Hamiltonian isotopy classes of Lagrangian fillings as there are cluster seeds in the $D_n$-type cluster algebra. This will be a consequence of a stronger geometric result, weave realizability in $D-$type, which we discuss below.

\begin{center}
		\begin{figure}[h!] 
			\centering
			\includegraphics[scale=0.8]{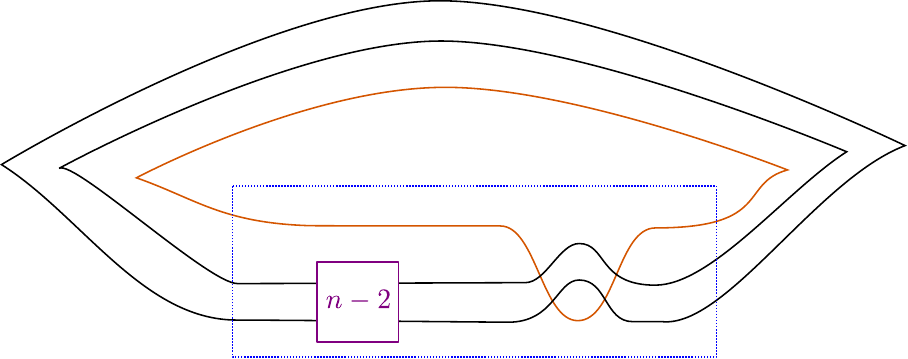}\caption{The front projection of $\lambda(D_n)\subseteq (\mathbb{S}^3,\xi_{\st})$. The box labelled with an $n-2$ represents $n-2$ positive crossings given by $\sigma_1^{n-2}.$ When $n$ is even, $\lambda(D_n)$ has 3-components, while when $n$ is odd, $\lambda(D_n)$ only has 2 components.} 
			\label{fig:DnLink}\end{figure}
	\end{center}

	By definition, the Legendrian link $\lambda(D_n)\subseteq(\mathbb{S}^3,\xi_{\st})$, $n\geq 4$ of $D_n$-type is the standard satellite of the Legendrian link defined by the front projection given by the 
	3-stranded positive braid $\sigma_1^{n-2}(\sigma_2\sigma_1^2\sigma_2)(\sigma_1\sigma_2)^3$, where $\sigma_1$ and $\sigma_2$ are the Artin generators for the 3-stranded braid group. Figure \ref{fig:DnLink} depicts a front diagram for $\lambda(D_n)$; note that the $(-1)$-framed closure of $\sigma_1^{n-2}(\sigma_2\sigma_1^2\sigma_2)(\sigma_1\sigma_2)^3$ is Legendrian isotopic to the rainbow closure of $\sigma_1^{n-2}(\sigma_2\sigma_1^2\sigma_2)$, the latter being depicted. The Legendrian link $\lambda(D_n)$ is also a max-tb representative of the smooth isotopy class of the link of the singularity $f(x,y)=y(x^2+y^{n-2})$. Since these are algebraic links, the max-tb representative given above is unique -- e.g. \cite[Proposition 2.2]{CasalsLagSkel} -- and has at least one exact Lagrangian filling \cite{HaydenSabloff}.

\begin{center}\begin{figure}[h!]{ \includegraphics[width=\textwidth]{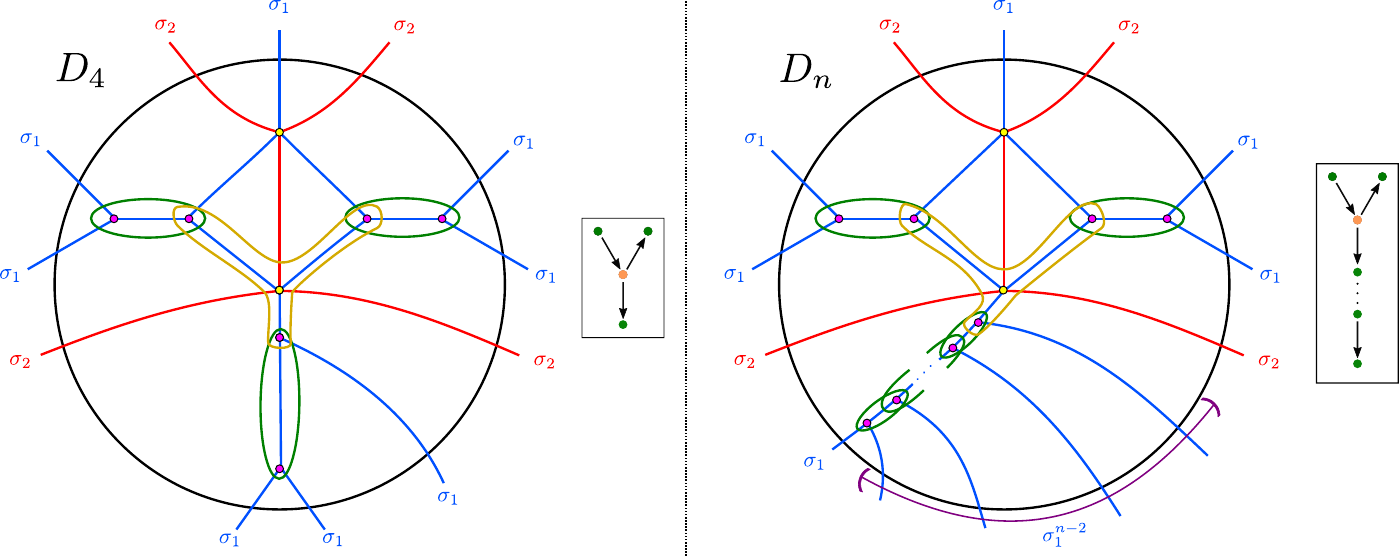}}\caption{3-graphs $\Gamma_0(D_4)$ (left) and $\Gamma_0(D_n)$ (right), pictured with their associated intersection quivers $Q(\Gamma_0(D_4), \{\gamma_i^{(0)}\})$ and $Q(\Gamma_0(D_n), \{\gamma_i^{(0)}\})$. The basis $\{\gamma_i^{(0)}\}$ for $H_1(\Lambda(\Gamma_0(D_4));\Z)$ is depicted by the dark green and orange and cycles drawn in the graph. Note that the quivers correspond to the $D_4$ and $D_n$ Dynkin diagrams, usually depicted rotated $90^\circ$ counterclockwise. 
			}
			\label{fig:D4Graph}\end{figure}
	\end{center}

	The $N$-graph calculus developed by Casals and Zaslow in \cite{CasalsZaslow} allows us to associate an exact Lagrangian filling of a ($-1$)-framed closure of a positive braid to a pair of trivalent planar graphs satisfying certain properties. See Figure \ref{fig:D4Graph} (left) for an example of a particular 3-graph, denoted by $\Gamma_0(D_4)$, associated to the Legendrian link $\lambda(D_4)$.\footnote{We use $\lambda(D_4)$, i.e. $n=4$, as a first example because $n=3$ would correspond to $\lambda(A_3)$, which has been studied previously \cite{EHK, Pan-fillings}. The study of $\lambda(D_4)$ is also the first instance where we require the machinery of 3-graphs rather than 2-graphs.} 
	In Section 3, we will show that the 3-graph $\Gamma_0(D_4)$ generalizes to a family of 3-graphs $\Gamma_0(D_n)$, depicted in Figure $\ref{fig:D4Graph}$ (right) for any $n\geq 3.$ 
	In a nutshell, a 3-fold branched cover of $\mathbb{D}^2$, simply branched at the trivalent vertices of these 3-graphs, yields an exact Lagrangian surface in $(T^*\mathbb{D}^2,\omega_{\st})$, whose Legendrian lift is a Legendrian weave. One of the distinct advantages of the 3-graph calculus is that it combinatorializes an operation, known as Lagrangian disk surgery \cite{Polterovich_Surgery,Yau}
	that modifies the weave in such a way as to yield additional -- non-Hamiltonian isotopic -- exact Lagrangian fillings of the link. 
	
	If we consider a 3-graph $\Gamma$ and a basis $\{\gamma_i\}$ for the first homology of the weave $\Lambda(\Gamma)$, $i\in[1,b_1(\Lambda(\Gamma))]$, we can define a quiver $Q(\Gamma, \{\gamma_i\})$ whose adjacency matrix is given by the intersection form in $H_1(\Lambda(\Gamma))$. Quivers come equipped with a involutive operation, known as quiver mutation, that produces new quivers; see subsection \ref{quivers} below or \cite{FWZ1} for more on quivers.  A key result of \cite{CasalsZaslow} tells us that Legendrian mutation of the weave induces a quiver mutation of the intersection quiver. Quivers related by a sequence of mutations are said to be mutation equivalent, and the quivers that are of finite mutation type (i.e. the set of mutation equivalent quivers is finite) have an ADE classification \cite{FWZ2}. This classification parallels the naming convention for the $D_n$ links described above: the intersection quiver associated to $\lambda(D_n)$ is a quiver in the mutation class of the $D_n$-Dynkin diagram (the latter endowed with an appropriate orientation). See Figure \ref{fig:D4Graph} for examples of $D_4$ and a $D_n$ quivers. For our 3-graph $\Gamma_0(D_n)$, $n\geq 3$, we will give an explicit basis $\{\gamma_i^{(0)}\}=\{\gamma_1^{(0)}, \ldots, \gamma_n^{(0)}\}$,  for $H_1(\Lambda(\Gamma_0(D_n)),\Z)$, whose intersection quiver $Q(\Gamma_0(D_n), \{\gamma_i^{(0)}\})$ is the standard $D_n$-Dynkin diagram. 
	
	Let us introduce the following notion in this manuscript. By definition, a sequence of quiver mutations for $Q(\Gamma_0(D_n), \{\gamma_i^{(0)}\})$ is said to be {\it weave realizable} if each quiver mutation in the sequence can be realized as a Legendrian weave mutation for a 3-graph. Our main result is the following theorem:


	\begin{theorem}\label{main}{
			Any sequence of quiver mutations of $Q(\Gamma_0(D_n), \gamma_1^{(0)}, \ldots, \gamma_n^{(0)}\})$ is weave realizable.}
	\end{theorem}
	
	In other words, Theorem \ref{main} states that in $D$-type, any algebraic quiver mutation can actually be realized geometrically by a Legendrian weave mutation. Weave realizability is of interest because it measures the difference between algebraic invariants -- e.g. the cluster structure in the moduli of sheaves -- and geometric objects, in this case Hamiltonian isotopy classes of exact Lagrangian fillings. In general if any sequence of quiver mutations were weave realizable, we would know that each cluster is inhabited by at least one embedded exact Lagrangian filling -- this general statement remains open for an arbitrary Legendrian link. For instance, any link with an associated quiver that is not of finite mutation type 
	satisfying the weave realizability property would admit infinitely many Lagrangian fillings, distinguished by their quivers.\footnote{This would be independent of the cluster structure defined by the microlocal monodromy functor, which we actually must use for $D$-type.} Note that weave realizability was shown for A-type in \cite{TreumannZaslow}, and beyond $A$- and $D$-types we currently do not know whether there are any other links satisfying the weave realizability property.

	

	We can further distinguish fillings by studying the cluster algebra structure on the moduli of microlocal rank-1 sheaves $\mathcal{C}(\Gamma)$ of a weave $\Lambda(\Gamma)$, e.g. see \cite{CasalsZaslow}. Specifically, sheaf quantization of each exact Lagrangian filling of $\lambda(D_n)$ induces a cluster chart on the coordinate ring of functions on $\mathcal{C}(\Gamma_0(D_n))$ via the microlocal monodromy functor, giving $\mathcal{C}(\Gamma_0(D_n))$ the structure of a cluster variety of $D_n$-type \cite{STZ_ConstrSheaves, STWZ}. Describing a single cluster chart in this cluster variety requires the data of the quiver associated to the weave, and the microlocal monodromy around each 1-cycle of the weave. Crucially, applying the Legendrian mutation operation to the weave induces a cluster transformation on the cluster chart, and the specific cluster chart defined by a Lagrangian fillings is a Hamiltonian isotopy invariant. Therefore, Theorem \ref{main} has the following consequence.

	\begin{corollary} \label{cor} 
		Every cluster chart of the moduli of microlocal rank-$1$ sheaves $\mathcal{C}(\Gamma_0(D_n))$ is induced by at least one embedded exact Lagrangian filling of $\lambda(D_n)\sse(\S^3,\xi_{\st})$.  In particular, there exist at least $(3n-2)C_{n-1}$ exact Lagrangian fillings of the link $\lambda(D_n)$ up to Hamiltonian isotopy, where $C_{n}$ denotes the $n$th Catalan number.
	\end{corollary}
	
	Moreover, weave realizability implies a slightly stronger result. Specifically, we can consider the {\it filling exchange graph} associated to a link of $D_n$-type, where the vertices are Hamiltonian isotopy classes of embedded exact Lagrangians, and two vertices are connected by an edge if the two fillings are related by a Lagrangian disk surgery. Then weave realizability implies that the filling exchange graph contains a subgraph 
	isomorphic to the cluster exchange graph for the cluster algebra of $D_n$-type.
	
	\begin{remark}
		As of yet, we have no way of determining whether our method produces all possible exact Lagrangian fillings of a type $D_n$-link. This question remains open for $A$-type Legendrian links as well. In fact, the only known knot for which we have a complete nonempty classification of Lagrangian fillings is the Legendrian unknot, which has a unique filling \cite{EliashbergPolterovich96}. 
		\hfill$\Box$
	\end{remark}

	In summary, our method for constructing exact Lagrangian fillings will be to represent them using the planar diagrammatic calculus of N-graphs developed in \cite{CasalsZaslow}. This diagrammatic calculus includes a mutation operation on the diagrams that yields additional fillings. We distinguish the resulting fillings up to Hamiltonian using a sheaf-theoretic invariant. From this data, we extract a cluster algebra structure and show that every mutation of the quiver associated to the cluster can be realized by applying our Legendrian mutation operation to the 3-graph, thus proving that there are at least as many distinct fillings as distinct cluster seeds of $D_n$-type. The main theorem will be proven in Section 3 after giving the necessary preliminaries in Section 2.   
	
	\subsection*{Acknowledgments}  Many thanks to Roger Casals for his support and encouragement throughout this project. Thanks also to Youngjin Bae and Eric Zaslow for helpful conversations, and to the anonymous referee for insightful comments.

	\subsection*{Added in proof}
	While writing this manuscript, we learned that recent independent work by Byung Hee An, Youngjin Bae, and Eunjeong Lee also produces at least as many exact Lagrangian fillings as cluster seeds for links of $ADE$ type \cite{ABL2021}, providing an alternative proof to Corollary 1. From our understanding, they use an inductive argument that relies on the combinatorial properties of the finite type generalized associahedron. Specifically, they leverage the fact that the Coxeter transformation in finite type is transitive if starting with a particular set of vertices by finding a weave pattern 
	that realizes Coxeter mutations. 
	While their initial 3-graph $\mathcal{G}(1, 1, n)$ is the same as our $\Gamma_0(D_n),$ their method of computing a weave associated to an arbitrary sequence of quiver mutations requires concatenating some number of concordances corresponding to the Coxeter mutation before mutating. As a result, a 3-graph arising from a sequence of quiver mutations $\mu_1\dots \mu_i$ computed using this method is not explicitly shown to be related to a 3-graph arising from a sequence of quiver mutations $\mu_1\dots \mu_i \mu_{i+1}$ by a single Legendrian mutation of the weave. In contrast, in our approach we are able to relate each 3-graph arising from a sequence of quiver mutations to the next by a single Legendrian mutation and a specific set of Legendrian Reidemeister moves. While both this manuscript and \cite{ABL2021} use the framework of $N$-graphs to approach the problem of enumerating exact Lagrangian fillings, the proofs are different, independent, and our approach is able to give an explicit construction for realizing any sequence of quiver mutations via an explicit sequence of mutations in the 3-graph.\hfill$\Box$


	\section{Preliminaries}

	In this section we introduce the necessary ingredients required for the proof of Theorem \ref{main} and Corollary \ref{cor}. We first discuss the contact topology needed to understand weaves and their homology. We then discuss the sheaf-theoretic material related to distinguishing fillings via cluster algebraic methods.
	
	
	\subsection{Contact Topology and Exact Lagrangian Fillings}
	
	A contact structure $\xi$ on $\R^3$ is a 2-plane field given locally as the kernel of a 1-form $\alpha\in \Omega^1(\R^3)$ satisfying $\alpha\wedge d\alpha\neq 0$. The standard contact structure on $(\R^3, \xi_{\st})$ is given by the kernel of $\alpha=dz - ydx$.  A Legendrian link
	$\lambda$ in $(\R^3, \xi)$ 
	is an embedding of a disjoint union of copies of $\mathbb{S}^1$ that is always tangent to $\xi$. By definition, the contact 3-sphere $(\S^3, \xi_{\st})$ is the one point compactification of $(\R^3,\xi_{\st})$ . Since a link in $\S^3$ can always be assumed to avoid a point, we will equivalently be considering Legendrian links in $(\R^3, \xi_{\st})$ and $(\S^3, \xi_{\st}).$ By definition, the symplectization of  $(\R^3, \xi_{\st})$ is given by $(\R^3\times \R_t, d(e^t \alpha))$. 
	
	Given two Legendrian links $\lambda_+$ and $\lambda_-$ in $(\R^3,\xi)$, an exact Lagrangian cobordism $\Sigma$ from $\lambda_-$ to $\lambda_+$ is an embedded compact orientable surface in the symplectization $(\R^3\times \R_t, d(e^t\alpha))$ such that for some $T>0$

	\begin{itemize}
		\item $\Sigma \cap \left(\R^3\times [T,\infty)\right)=\lambda_+\times [T, \infty)$
		\item $\Sigma \cap \left(\R^3 \times (-\infty, -T)\right)=\lambda_-\times (-\infty, -T]$
		\item $\Sigma$ is an exact Lagrangian, i.e. $e^t\alpha=df$ for some function $f:\Sigma\to \R.$
	\end{itemize}
	The asymptotic behavior of $\Sigma$, as specified by the first two conditions, ensures that we can concatenate Lagrangian cobordisms.
	By definition, an exact Lagrangian filling of $\lambda_+$ is an exact Lagrangian cobordism from $\emptyset$ to $\lambda_+$.
	

	
	We can also consider the Legendrian lift of an exact Lagrangian in the contactization $(\R_s\times \R^4,\ker\{ds-d(e^t\alpha)\})$ of $(\R^4, d(e^t\alpha))$. Note that there exists a contactomorphism between $(\R_s\times \R^4,\ker\{ds-d(e^t\alpha)\})$ and the standard contact Darboux structure $(\R^5,\xi_{\st})$, where $\xi_{\st}=\ker\{dz-y_1dx_1-y_2dx_2\}$. We will often work with the Legendrian front projection $(\R^5,\xi_{\st})\longrightarrow\R^3_{x_1,x_2,z}$ for the latter. This will be a useful perspective for us, as it allows us to construct Lagrangian fillings by studying (wave)fronts in $\R^3=\R^3_{x_1,x_2,z}$ of Legendrian surfaces in $(\R^5,\xi_{\st})$, and then projecting down to the standard symplectic Darboux chart $\R^4=\R^4_{x_1,y_1,x_2,y_2}$. In this setting, the exact Lagrangian surface is embedded in $\R^4$ if and only if its Legendrian lift has no Reeb chords. The construction will be performed through the combinatorics of $N$-graphs, as we now explain.

	
	\subsection{3-graphs and Weaves}

	In this subsection, we discuss the diagrammatic method of constructing and manipulating exact Lagrangian fillings of links arising as the ($-1$)-framed closures of positive braids via the calculus of $N$-graphs. For this manuscript, it will suffice to take $N=3$.
	\begin{definition}
		A 3-graph is a pair of 
		embedded planar trivalent graphs $B,R\subseteq \D^2$ such that at any vertex $v\in B\cap R$ the six edges belonging to $B$ and $R$ incident to $v$ alternate. \hfill $\Box$ 
	\end{definition}
	Equivalently, a 3-graph is an edge-bicolored graph with monochromatic trivalent vertices and interlacing hexavalent vertices. $\Gamma_0(D_4),$ depicted in Figure \ref{fig:D4Graph} (left) contains two hexavalent vertices displaying the alternating behavior described in the definition.

	\begin{remark}
		\cite{CasalsZaslow} gives a general framework for working with N-graphs, where $N-1$ is the number of embedded planar trivalent graphs. This allows for the study of fillings of Legendrian links associated to $N$-stranded positive braids. This can also be generalized to consider N-graphs in a surface other than $\D^2$. In our case, the family of links $\la(D_n)$ can be expressed as a family of 3-stranded braids, hence our choice to restrict $N$ to 3 in $\D^2$. \hfill $\Box$
	\end{remark}

	Given a 3-graph $\Gamma \subseteq \mathbb{D}^2,$ we describe how to associate a Legendrian surface $\Lambda(\Gamma)\subseteq (\R^5, \xi_{\st})$. To do so, we first describe certain singularities of $\Lambda(\Gamma)$ that arise under the Legendrian front projection $\pi:(\R^5, \xi_{\st})\to (\R^3,\xi_{\st})$. In general, such singularities are known as Legendrian singularities or singularities of fronts. See 
	\cite{ArnoldSing} for a classification of such singularities. The three singularities we will be interested in are the $A_1^2$, $A_1^3$ and $D_4^-$ singularities, pictured in Figure \ref{fig:SingWave} below.  

	\begin{center}
		\begin{figure}[h!]{ \includegraphics[width=.9\textwidth]{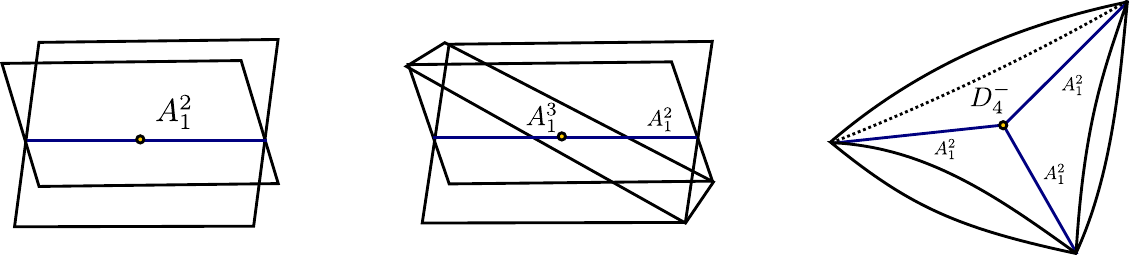}}\caption{$A_1^2$ (left), $A_1^3$ (center), and $D_4^-$ (right) singularities represented in the 3-graph by an edge, hexavalent vertex, and trivalent vertex, respectively.}
			\label{fig:SingWave}\end{figure}
	\end{center}

	Before we describe our Legendrian surfaces, we must first discuss the ambient contact structure that they live in. For $\Gamma\subseteq \D^2$ we will take $\Lambda(\Gamma)$ to live in the first jet space $(J^1\D^2,\xi_{\st})=(T^*\D^2\times \R_z, \ker(dz-\theta_{\st}))$, where $\theta_{\st}$ is the standard Liouville form on the cotangent bundle $T^*\D^2$. We can view $J^1\D^2$ as a certain local model for a contact structure, in the following way. If we take $(Y,\xi)$ to be a contact 5-manifold, then by the Weinstein neighborhood theorem, any Legendrian embedding $i:\D^2\to (Y,\xi)$ extends to an embedding from $(J^1\D^2,\xi_{\st})$ to a small open neighborhood of $i(\D^2)$ with contact structure given by the restriction of $\xi$ to that neighborhood. In particular, a Legendrian embedding of $i:\mathbb{S}^1 \to \mathbb{S}^3$ gives rise to a contact embedding $\tilde{i}:J^1\mathbb{S}^1\longrightarrow\mbox{Op}(i(\S^1))$ into some open neighborhood $\mbox{Op}(i(\S^1))\subseteq\mathbb{S}^3$. Of particular note in our case is that, under a Legendrian embedding $\D^2\subseteq(\R^5,\xi_{\st})$, a Legendrian link $\lambda$ in $J^1\dd\D^2$ is mapped to a Legendrian link in the contact boundary $(\S^3,\xi_{\st})$ of the symplectic $(\R^4,\omega_{\st}=d\theta_{\st})$ given as the co-domain of the Lagrangian projection $(\R^5,\xi_{\st})\lr(\R^4,\omega_st)$. See \cite{Satellite2} for a description of this Legendrian satellite operation. 
	
	To construct a Legendrian weave $\Lambda(\Gamma)\subseteq (J^1\D^2,\xi_{\st})$ from a 3-graph $\Gamma$, we glue together the local germs of singularities according to the edges of $\Gamma$. First, consider three horizontal wavefronts $\D^2\times \{1\}\sqcup \D^2\times \{2\}\sqcup \D^2\times \{3\}\subseteq \D^2\times \R$ and a 3-graph $\Gamma\subseteq \D^2\times \{0\}$. We construct the associated Legendrian weave $\Lambda(\Gamma)$ as follows.
	\begin{itemize}
		\item Above each blue (resp. red) edge, insert an $A_1^2$ crossing between the $\D^2\times\{1\}$ and $\D^2\times \{2\}$ sheets (resp $\D^2\times\{2\}$ and $\D^2\times\{3\}$ sheets) so that the projection of the $A_1^2$ singular locus under $\pi:\D^2\times \R \to \D^2\times \{0\}$ agrees with the blue (resp. red) edge.   
		\item At each blue (resp. red) trivalent vertex $v$, insert a $D_4^-$ singularity between the sheets $\D^2\times \{1\}$ and $\D^2\times\{2\}$ (resp. $\D^2\times\{2\}$ and $\D^2\times\{3\}$) in such a way that the projection of the $D_4^-$ singular locus agrees with $v$ and the projection of the $A_2^1$ crossings agree with the edges incident to $v$.
		\item At each hexavalent vertex $v$, insert an $A_1^3$ singularity along the three sheets in such a way that the origin of the $A_1^3$ singular locus agrees with $v$ and the $A_1^2$ crossings agree with the edges incident to $v$.
	\end{itemize} 
	\begin{center}
		\begin{figure}[h!]{ \includegraphics[width=\textwidth]{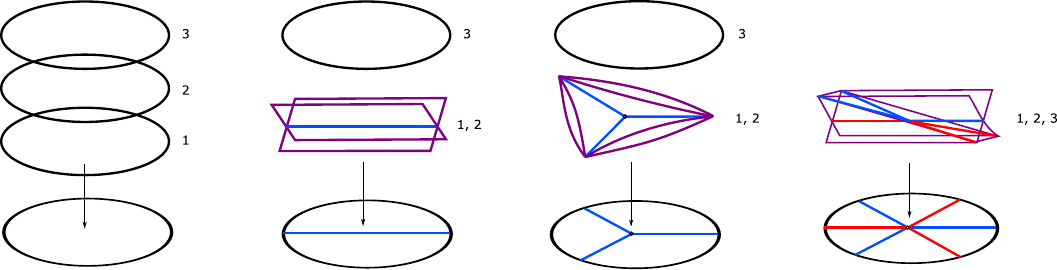}}\caption{The weaving of the singularities pictured in Figure \ref{fig:SingWave} along the edges of the N-graph. Gluing these local pictures together according to the 3-graph $\Gamma$ yields the weave $\Lambda(\Gamma)$.}
			\label{fig:Weaving}\end{figure}
	\end{center}
	If we take an open cover $\{U_i\}_{i=1}^m$ of $\D^2\times \{0\}$ by open disks, refined so that any disk contains at most one of these three features, we can glue together the resulting fronts according to the intersection of edges along the boundary of our disks. Specifically, if $U_i\cap U_j$ is nonempty, then we define $\Sigma(U_1\cup U_2)$ to be the wavefront resulting from considering the union of wavefronts $\Sigma(U_1)\cup\Sigma(U_j)$ in $(U_1\cup U_2)\times \R$. We define the Legendrian weave $\Lambda(\Gamma)$ as the Legendrian surface contained in $(J^1\D^2, \xi_{\st})$ with wavefront $\Sigma(\Gamma)=\Sigma(\cup_{i=1}^m U_i)$ given by gluing the local wavefronts of singularities together according to the 3-graph $\Gamma$ \cite[Section 2.3]{CasalsZaslow}.
	
	The smooth topology of a Legendrian weave $\Lambda(\Gamma)$ is given 
	as a 3-fold branched cover over $\D^2$ with simple branched points corresponding to each of the trivalent vertices of $\Gamma$. 
	The genus of $\Lambda(\Gamma)$ is then computed using the Riemann-Hurwitz formula: 
	
	$$
	g(\Lambda(\Gamma))=\frac{1}{2}(v(\Gamma)+2-3\chi(\D^2)-|\partial\Lambda(\Gamma)|)$$
	where $v(\Gamma)$ is the number of trivalent vertices of $\Gamma$ and $|\partial\Lambda(\Gamma)|$ denotes the number of boundary components of $\Gamma$. 
	
	\begin{ex}
		If we apply this formula to the 3-graph $\Gamma_0(D_4)$, pictured in Figure \ref{fig:D4Graph}, we have $6$ trivalent vertices and 3 link components, so the genus is computed as 
		$g(\Lambda(\Gamma_0(D_4)))=\frac{1}{2}(6 +2 - 3 -3 )=1.$
		
		For $\Gamma_0(D_n)$, we have three boundary components for even $n$ and two boundary components for odd n. The number of trivalent vertices is $n+2$, so the genus $g(\Lambda(\Gamma_0(D_n))$ is $\lfloor\frac{n-1}{2}\rfloor$, assuming $n\geq 2$.

		This computation tells us that $\Lambda(\Gamma_0(D_4))$ is smoothly a 3-punctured torus bounding the link $\lambda(D_4).$ Therefore, we can give a basis for $H_1(\Lambda(\Gamma_0(D_4));\Z)$ in terms of the four cycles pictured in Figure \ref{fig:D4Graph}. \end{ex}
	For $\Gamma_0(D_n)$, the corresponding weave $\Lambda(\Gamma_0(D_n))$ will be smoothly a genus $\lfloor\frac{n-1}{2}\rfloor$ surface with a basis of $H_1(\Lambda(\Gamma);\Z)$ given by $n$ cycles. Our computation of the genus in the example above agrees with a theorem of Chantraine \cite{Chantraine10} specifying the relationship between the Thurston-Bennequin invariant of $\la(D_n)$ and the genus of any exact Lagrangian filling $\Lambda$ of $\la(D_n)$. In particular, $tb(\la(D_n))=n-1$ and therefore the Euler characteristic of $\Lambda$ is $3-n$ when $n$ is odd and $4-n$ when $n$ is even. Thus, we recover the genus $g(\Lambda)=\lfloor\frac{n-1}{2}\rfloor$ of any filling of $\la(D_n)$. 
	In the next section, we describe a general method for giving a basis $\{\gamma_i^{(0)}\}, i\in [1, n]$ of the first homology $H_1(\La(\Gamma_0(D_n)); \Z)\cong \Z^n$.
	\

	

	\subsection{Homology of Weaves}
\begin{center}
		\begin{figure}[h!]{ \includegraphics[width=.95\textwidth]{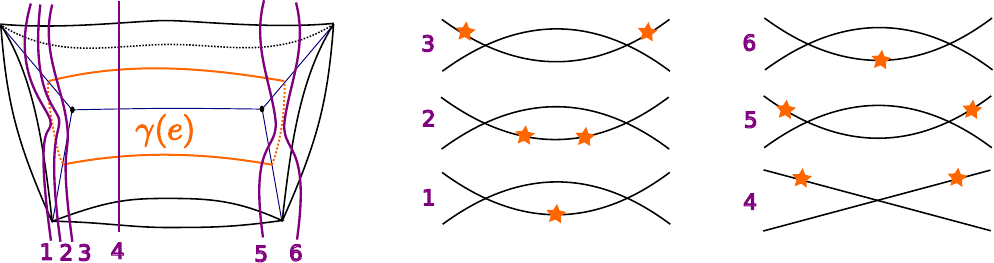}}\caption{A short {\sf I}-cycle $\gamma(e)$ for the edge $e\in G$ pictured in the wavefront $\Sigma(\Gamma)$ (left) and a vertical slicing of $\Sigma(\Gamma)$ (right).}
			\label{fig:I-cycle}\end{figure}
	\end{center}

\begin{center}\begin{figure}[h!]{ \includegraphics[width=.95\textwidth]{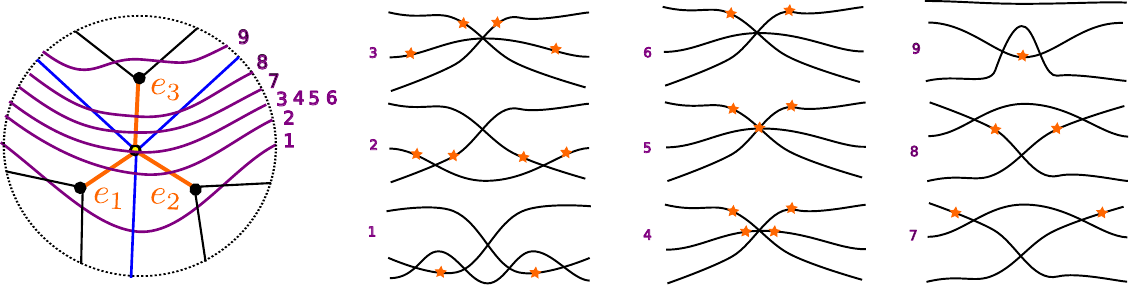}}\caption{A short {\sf Y}-cycle $\gamma(e)$ defined by the edges $e_1, e_2, e_3\in G$ pictured in the wavefront $\Sigma(\Gamma)$ (left) and a vertical slicing of $\Sigma(\Gamma)$ (right).}
			\label{fig:Y-cycle}\end{figure}
	\end{center}

\begin{center}
		\begin{figure}[h!]{ \includegraphics[width=.95\textwidth]{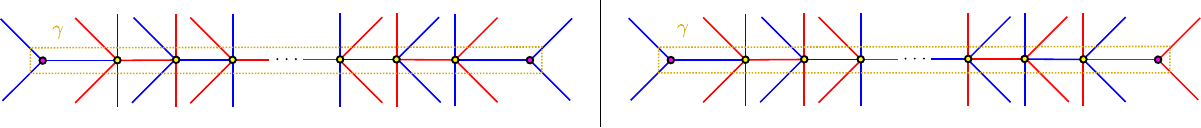}}\caption{A pair of long {\sf I}-cycles, both denoted by $\gamma$. The cycle on the left passes through an even number of hexavalent vertices, while the cycle on the right passes through an odd number. 
			}
			\label{fig:longI-cycle}\end{figure}
	\end{center}

	We require a description of the first homology $H_1(\Lambda(\Gamma));\Z)$ in order to apply the mutation operation to a 3-graph $\Gamma$. We first consider an edge connecting two trivalent vertices. Closely examining the sheets of our surface, we can see that each such edge corresponds to a 1-cycle, as pictured in Figure \ref{fig:I-cycle} (left). We refer to such a 1-cycle as a short {\sf I}-cycle. Similarly, any three edges of the same color that connect a single hexavalent vertex to three trivalent vertices correspond to a 1-cycle, as pictured in \ref{fig:Y-cycle} (left). We refer to such a 1-cycle as a short {\sf Y}-cycle. See figures \ref{fig:I-cycle} (right) and \ref{fig:Y-cycle} (right) for a diagram of these 1-cycles in the wavefront $\Sigma(\Gamma)$. We can also consider a sequence of edges starting and ending at trivalent vertices and passing directly through any number of hexavalent vertices, as pictured in Figure \ref{fig:longI-cycle}. Such a cycle is referred to as a long {\sf I}-cycle. Finally, we can combine any number of {\sf I}-cycles and short {\sf Y}-cycles to describe an arbitrary 1-cycle as a tree with leaves on trivalent vertices and edges passing directly through hexavalent vertices.
	
\begin{center}\begin{figure}[h!]{ \includegraphics[width=.5\textwidth]{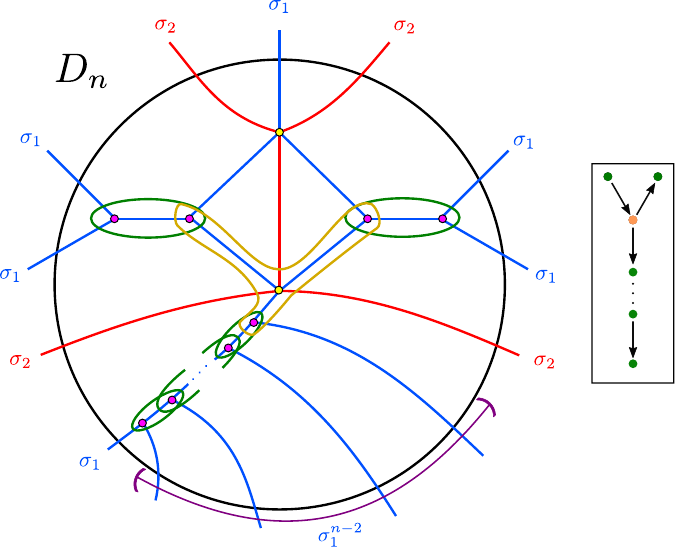}}\caption{The 3-graph $\Gamma_0(D_n)$ and its associated intersection quiver. The basis $\{\gamma_i^{(0)}\}$ of $H_1(\Lambda(\Gamma_0(D_n));\Z)$ is given by the orange {\sf Y}-cycle, the green {\sf I}-cycles, and the $n-3$ {\sf I}-cycles not pictured.}
			\label{fig:DnGraph}\end{figure}
	\end{center}

	In the proof of our main result, we will generally give a basis for $H_1(\Lambda(\Gamma);\Z)$ in terms of short {\sf I}-cycles and short {\sf Y}-cycles. Indeed, Figure \ref{fig:DnGraph} gives a basis of $H_1(\Lambda(\Gamma_0(D_n));\Z)$ consisting of $n-1$ short {\sf I}-cycles and a single {\sf Y}-cycle. 

\begin{center}\begin{figure}[h!]{ \includegraphics[width=.4\textwidth]{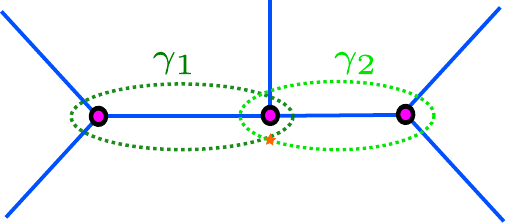}}\caption{Intersection of two cycles, $\gamma_1$ and $\gamma_2$. The intersection point is indicated by an orange star. If we orient both cycles counterclockwise, then we will set $\langle \gamma_1,\gamma_2\rangle=-1$ as our convention.} 
			\label{fig:intersect}\end{figure}
	\end{center}

	The	intersection form $\langle \cdot, \cdot \rangle$ on $H_1(\Lambda(\Gamma))$ plays a key role in distinguishing our Legendrian weaves. If we consider a pair of 1-cycles $\gamma_1, \gamma_2\in H_1(\Lambda(\Gamma))$ with nonempty geometric intersection in $\Gamma$, as pictured in Figure \ref{fig:intersect}, we can see that the intersection of their projection onto the 3-graph differs from the intersection in $\Lambda(\Gamma).$ Specifically, we can carefully examine the sheets that the 1-cycles cross in order to see that $\gamma_1$ and $\gamma_2$  intersect only in a single point of $\Lambda(\Gamma)$. If we fix an orientation on $\gamma_1$ and $\gamma_2,$ then we can assign a sign to this intersection based on the convention given in Figure \ref{fig:intersect}. We refer to the signed count of the intersection of $\gamma_1$ and $\gamma_2$ as their algebraic intersection and denote it by $\langle \gamma_1, \gamma_2\rangle.$ For the remainder of this manuscript, we will fix a counterclockwise orientation for all of our cycles and adopt the convention that any two cycles $\gamma_1$ and $\gamma_2$, intersecting at a trivalent vertex as in Figure \ref{fig:intersect} have algebraic intersection $\langle \gamma_1, \gamma_2\rangle =-1.$
	
	{\bf Notation:} For the sake of visual clarity, we will represent an element of $H_1(\Lambda(\Gamma);\Z)$ by a colored edge for the remainder of this manuscript. This also ensures that the geometric intersection more accurately reflects the algebraic intersection. The original coloring of the blue or red edges can be readily obtained by examining $\Gamma$ and its trivalent vertices.\hfill$\Box$ 
	
\begin{center}\begin{figure}[h!]{ \includegraphics[width=\textwidth]{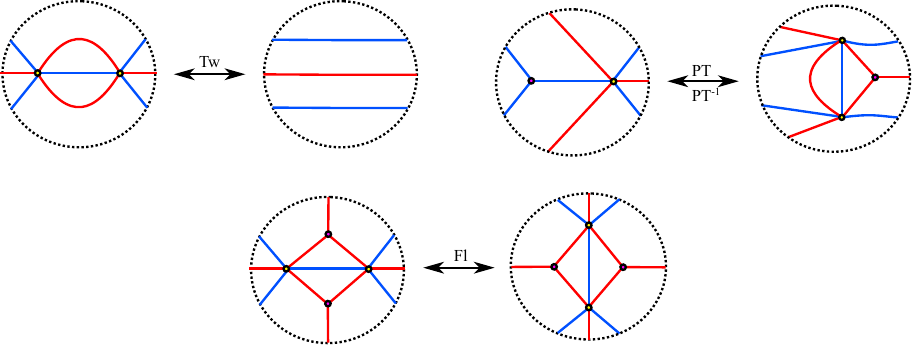}}\caption{Legendrian Surface Reidemeister moves for 3-graphs. From left to right, a candy twist, a push-through, and a flop, denoted by Tw, PT, and Fl respectively.} 
			\label{fig:Moves}\end{figure}
	\end{center}

	In our correspondence between 3-graphs and 
	weaves, we must consider how a Legendrian isotopy of the weave $\Lambda(\Gamma)$ affects the 3-graph $\Gamma$ and its homology basis. We can restrict our attention to certain isotopies, referred to as Legendrian Surface Reidemeister moves. These moves create specific changes in the Legendrian front $\Sigma(\Gamma)$, known as perestroikas or Reidemeister moves \cite{ArnoldSing}.  
	From \cite{CasalsZaslow}, we have the following theorem relating perestroikas of fronts to the corresponding 3-graphs. 
	\begin{theorem}[\cite{CasalsZaslow}, Theorem 4.2]
		Let $\Gamma$ and $\Gamma'$ be two 3-graphs related by one of the moves shown in Figure \ref{fig:Moves}. Then the associated weaves $\Lambda(\Gamma)$ and $\Lambda(\Gamma')$ are Legendrian isotopic relative to their boundaries. \hfill $\Box$
	\end{theorem}

\begin{center}\begin{figure}[h!]{ \includegraphics[width=\textwidth]{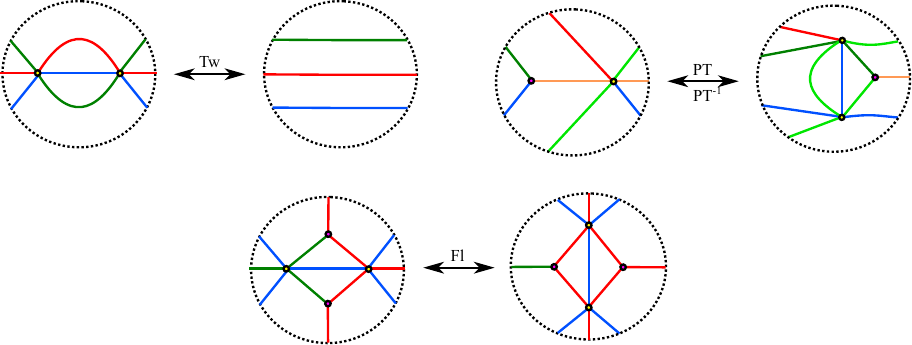}}\caption{Behavior of certain homology cycles under Legendrian Surface Reidemeister moves. 
			}\label{fig:ReidemeisterCycles}\end{figure}
	\end{center}

	See Figure \ref{fig:ReidemeisterCycles} for a description of the behavior of elements of $H_1(\Lambda(\Gamma); \Z)$ under these Legendrian Surface Reidemeister moves. In the pair of 3-graphs in Figure \ref{fig:ReidemeisterCycles} (center), we have denoted a push-through by PT or PT$^{-1}$ depending on whether we go from left to right or right to left.
	This helps us to specify the simplifications we make in the figures in the proof of Theorem \ref{main}, as this move is not as readily apparent as the other two. We will refer to the PT$^{-1}$ move as a reverse push-through. Note that an application of this move eliminates the geometric intersection between the light green and dark green cycles in Figure \ref{fig:ReidemeisterCycles}. 
	
	\begin{remark}
		It is also possible to verify the computations in Figure \ref{fig:ReidemeisterCycles} by examining the relative homology of a cycle. Specifically, if we have a basis of the relative homology $H_1(\Lambda(\Gamma), \partial\Lambda(\Gamma) ;\Z)$, then the intersection form on that basis allows us to determine a given cycle by Poincar\'e-Lefschetz duality. $\hfill \Box$
	\end{remark}
	
	
	\subsection{Mutations of 3-graphs}
	We complete our discussion of general 3-graphs with a description of Legendrian mutation, which we will use to generate distinct exact Lagrangian fillings. Given a Legendrian weave $\Lambda(\Gamma)$ and a 1-cycle $\gamma\in H_1(\Lambda(\Gamma);\Z)$, the Legendrian mutation $\mu_\gamma(\Lambda(\Gamma))$ outputs a 3-graph and a corresponding Legendrian weave smoothly isotopic to $\Lambda(\Gamma)$ but whose Lagrangian projection is generally not Hamiltonian isotopic to that of $\Lambda(\Gamma)$.
	\begin{definition}
		Two Legendrian surfaces $\Lambda_0, \Lambda_1 \subseteq (\R^5, \xi_{\st})$ with equal boundary $\partial \Lambda_0=\partial \Lambda_1$, are mutation-equivalent if and only if there exists a compactly supported Legendrian isotopy $\{\tilde{\Lambda}_t\}$ relative to the boundary, with $\tilde \Lambda_0=\Lambda_0$ and a Darboux ball $(B,\xi_{\st})$ such that 
		\begin{enumerate}
			\item[(i)] Outside the Darboux ball, we have $\tilde{\Lambda}_1|_{\R^5\backslash B}=\Lambda_1|_{\R^5\backslash B}$
			\item[(ii)] There exists a global front projection $\pi:\R^5\to \R^3$ such that the pair of fronts $\pi|_{B\cap \tilde{\Lambda_1}}$ and $\pi|_{B\cap \Lambda_1}$ coincides with the pair of fronts in Figure \ref{fig:LegendrianMutationFront} below.
		\end{enumerate} \hfill $\Box$
	\end{definition}

	\begin{center}\begin{figure}[h!]{ \includegraphics[width=.6\textwidth]{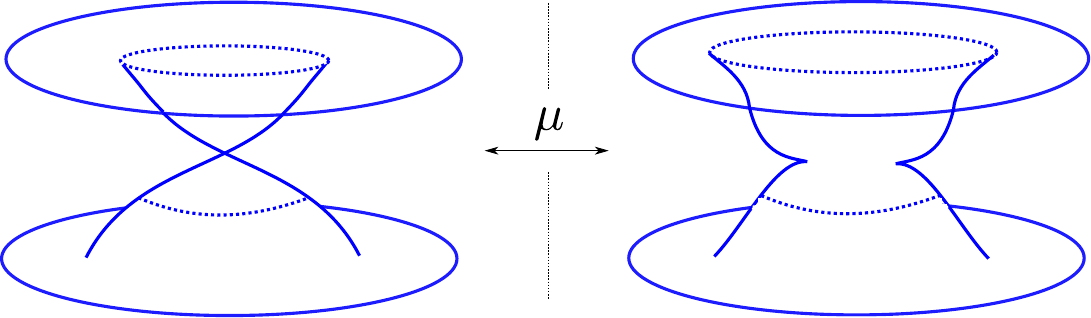}}\caption{Local fronts for two Legendrian cylinders non-Legendrian isotopic relative to their boundaries. 
			}\label{fig:LegendrianMutationFront}\end{figure}
	\end{center}
	
	We briefly note that these two fronts lift to non-Legendrian isotopic Legendrian cylinders in $(\R^5, \xi_{\st})$, relative to the boundary, and that the 1-cycle we input for our operation is precisely the 1-cycle defined by the cylinder corresponding to $\Lambda_0$.
	
\begin{center}\begin{figure}[h!]{ \includegraphics[width=.9\textwidth]{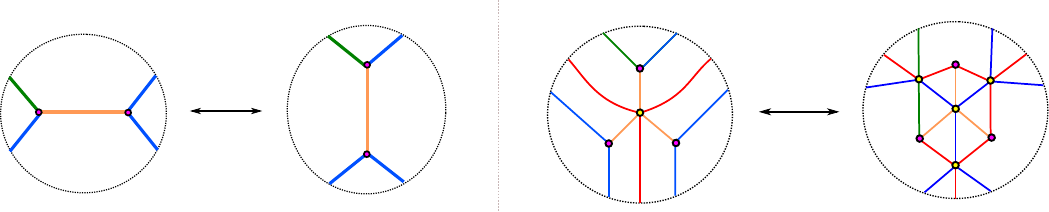}}\caption{Mutations of a 3-graph. The pair of 3-graphs on the left depicts mutation at the orange {\sf I}-cycle, while the pair of 3-graphs on the right depicts mutation at the orange {\sf Y}-cycle. In both cases, the dark green edge depicts the effect of mutation on any cycle intersecting the orange cycle.}\label{fig:Mutations}\end{figure}
	\end{center}

\begin{center}\begin{figure}[h!]{ \includegraphics[width=.8\textwidth]{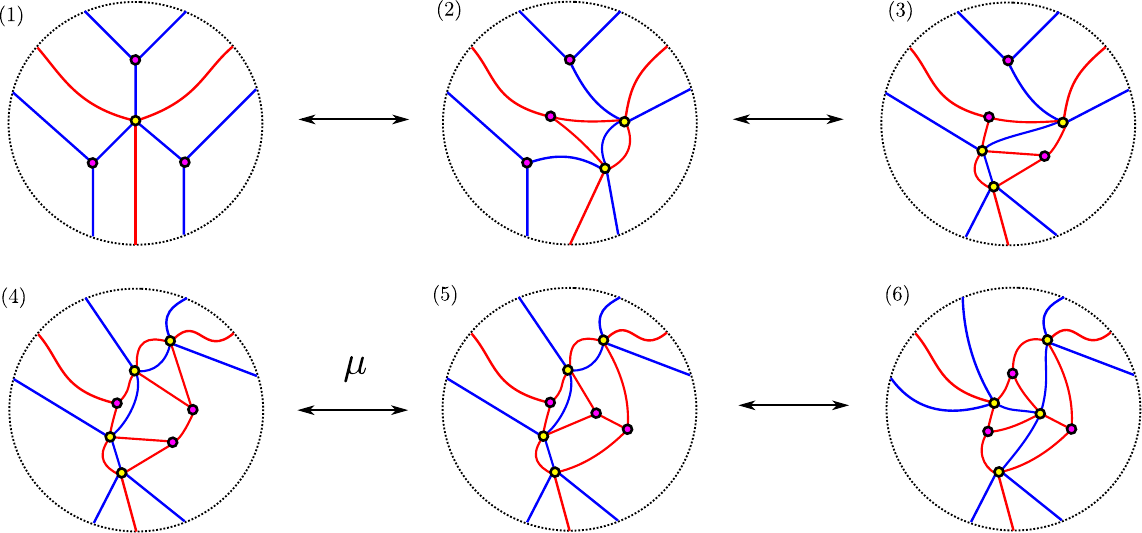}}\caption{Mutation at a short {\sf Y}-cycle given as a sequence of Legendrian Surface Reideister moves and mutation at a short {\sf I}-cycle. The {\sf Y}-cycle in the initial 3-graph is given by the three blue edges that each intersect the yellow vertex in the center.}\label{fig:Y-cycleMutation}\end{figure}
	\end{center}

	Combinatorially, we can describe mutation as certain manipulations of the edges of our graph. Figure \ref{fig:Mutations} (left) depicts mutation at a short {\sf I}-cycle, while Figure \ref{fig:Mutations} (right) depicts mutation at a short {\sf Y}-cycle. In the $N=2$ setting, we can identify 2-graphs with triangulations of an $n-$gon, in which case mutation at a short {\sf I}-cycle corresponds to a Whitehead move. 
	In the 3-graph setting, in order to describe mutation at a short {\sf Y}-cycle, we can first reduce the short {\sf Y}-cycle case to a short {\sf I}-cycle, as shown in Figure \ref{fig:Y-cycleMutation}, before applying our mutation.  See \cite[Section 4.9]{CasalsZaslow} for a more general description of mutation at long {\sf I}- and {\sf Y}-cycles in the 3-graph.

	The geometric operation above coincides with the combinatorial manipulation of the 3-graphs. Specifically, we have the following theorem.
	\begin{theorem}[\cite{CasalsZaslow}, Theorem 4.2.1]
		Given two 3-graphs, $\Gamma$ and $\Gamma'$ related by either of the combinatorial moves described in Figure \ref{fig:Mutations}, the corresponding Legendrian weaves $\Lambda(\Gamma)$ and $\Lambda(\Gamma')$ are mutation-equivalent relative to their boundary. \hfill $\Box$
	\end{theorem}
	

	\subsection{Lagrangian Fillings from Weaves}
	We now describe in more detail how an exact Lagrangian filling of a Legendrian link arises from a Legendrian weave. If we label all edges of $\Gamma\subseteq \D^2$ colored blue by $\sigma_1$ and all edges colored red by $\sigma_2$, then the points in the intersection $\Gamma\cap \partial \D^2$ give us a braid word in the Artin generators $\sigma_1$ and $\sigma_2$ of the 3-stranded braid group. We can then view the corresponding link $\beta$ as living in $(J^1 \mathbb{S}^1, \xi_{\st})$. If we consider our Legendrian weave $\Lambda(\Gamma)$ as an embedded Legendrian surface in $(\mathbb{R}^5,\xi_{\st})$, then according to our discussion above, it has boundary $\Lambda(\beta),$ where $\Lambda(\beta)$ is the Legendrian satellite of $\beta$ with companion knot given by the standard unknot. In our local contact model, the projection  $\pi:(J^1\D^2,\xi_{\st}) \to (T^*\D^2, \omega_st)$ gives an immersed exact Lagrangian surface with immersion points corresponding to Reeb chords of $\Lambda(\Gamma)$. If $\Lambda(\Gamma)$ has no Reeb chords, then $\pi$ is an embedding and $\Lambda(\Gamma)$ is an exact Lagrangian filling of $\Lambda(\beta).$ Since $(\mathbb{S}^3,\xi_{\st})$ minus a point is contactomorphic to $(\R^3, \xi_{\st})$, we have that an embedding of $\Lambda(\Gamma)$ into $(\mathbb{R}^5,\xi_{\st})$ gives an exact Lagrangian filling in $(\R^4,\xi_{\st})$ of $\Lambda(\beta)\subseteq (\R^3, \xi_{\st})$, as it can be assumed -- after a Legendrian isotopy -- to be disjoint from the point at infinity.
	
	\begin{remark}
		We study embedded -- rather than immersed -- Lagrangian fillings due to the existence of an $h$-principle for immersed Lagrangian fillings \cite[Theorem 16.3.2]{EliashbergMishachev02}. In particular, any pair of immersed exact Lagrangian fillings is connected by a one-parameter family of immersed exact Lagrangian fillings relative to the boundary. See also \cite{Gromov86}. 
	\end{remark}

	
	Our desire for embedded Lagrangians motivates the following definition.
	\begin{definition}
		A 3-graph $\Gamma\subseteq \D^2$ is free if the associated Legendrian front $\Sigma(\Gamma)$ can be woven with no Reeb chords. \hfill $\Box$
	\end{definition}
	
\begin{center}
		\begin{figure}[h!]{ \includegraphics[width=.8\textwidth]{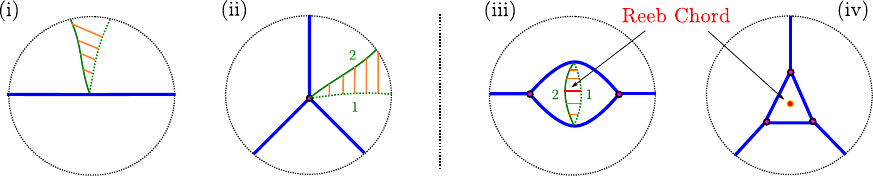}}\caption{2-graphs with woven with a choice of fronts illustrated by the green curves while the solid orange lines illustrate the difference in heights between sheets. A woven front for the pair of 2-graphs on the left can be chosen in such a way that the function giving the difference of heights between the two sheets of the front is 0 on $G$ and increasing towards the boundary. Critical points of the difference function correspond to Reeb chords, so the pair of 2-graphs on the left are free. However, any difference function for the pair of 2-graphs on the right must have at least one critical point inside the face.}\label{fig:Reebs}\end{figure}
	\end{center}

	In the $N=2$ setting, a 2-graph $\Gamma\subseteq \D^2$ is free if and only if $G$ has no bounded faces contained in the interior of $\D^2$. See figure \ref{fig:Reebs} for examples illustrating this characterization.   In the $N=3$ setting, there is no such simple characterization, but many 3-graphs can be determined to be free by direct inspection, as done in \cite[Section 7]{CasalsZaslow}. As an example, the 3-graph $\Gamma_0(D_n)$, depicted in Figure \ref{fig:DnGraph}, is a free 3-graph of $D_n$-type. This can be verified by taking a woven front for $\Lambda_0(D_n)$ such that the functions giving the difference of heights between the three sheets take the value 0 on $G$ and increase radially towards the boundary. Critical points of these difference functions correspond to Reeb chords. By construction, none of these difference functions have critical points, so $\Gamma_0(D_n)$ can be woven without Reeb chords and is a free 3-graph.

	Crucially, the mutation operation described above preserves the free property of a 3-graph.
	\begin{lemma}[\cite{CasalsZaslow}, lemma 7.4]
		Let $\Gamma\subseteq \D^2$ be a free 3-graph. Then the 3-graph $\mu(\Gamma)$ obtained by mutating according to any of the Legendrian mutation operations given above is also a free 3-graph. \hfill $\Box$
	\end{lemma}
	Therefore, starting with a free 3-graph and performing the Legendrian mutation operation gives us a method of creating additional embedded exact Lagrangian fillings. 
	
	At this stage, we have described the geometric and combinatorial ingredients needed for Theorem \ref{main}. The two subsequent subsections introduce the necessary algebraic invariants relating Legendrian weaves and 3-graphs to cluster algebras. These will be used to distinguish exact Lagrangian fillings.
	
	
	\subsection{Quivers from Weaves}\label{quivers}
	
	Before we describe the cluster algebra structure associated to a weave, we must first describe quivers and how they arise via the intersection form on $H_1(\Lambda(\Gamma);\Z).$ A quiver is a directed graph without loops or directed 2-cycles. In the Legendrian weave setting, the data of a quiver can be extracted from a given weave and a basis of its first homology. The intersection quiver is defined as follows: each basis element $\gamma_i\in H_1(\Lambda(\Gamma);\Z)$ defines a vertex $v_i$ in the quiver and we have $k$ arrows pointing from $v_j$ to $v_i$ if $\langle \gamma_i, \gamma_j\rangle=k$. We will only ever have $k$ either 0 or 1 for quivers arising from fillings of $\lambda(D_n)$. See Figure \ref{fig:D4Graph} (left) for an example of the quiver $Q(\Lambda(\Gamma_0(D_4)), \{\gamma_i^{(0)}\})$ defined by $\Lambda(\Gamma_0(D_4))$ and the indicated basis for $H_1(\Lambda(\Gamma_0(D_4);\Z)$.
	
	The combinatorial operation of quiver mutation at a vertex $v$ is defined as follows, e.g. see \cite{FWZ1}. First, for every pair of incoming edges and outgoing edges, we add an edge starting at the tail of the incoming edge and ending at the head of the outgoing edge. Next, we reverse the direction of all edges adjacent to $v$. Finally, we cancel any directed 2-cycles. If we started with the quiver $Q$, then we denote the quiver resulting from mutation at $v$ by $\mu_v(Q).$ See Figure \ref{fig:D4Mutation} (bottom) for an example. Under this operation, we can naturally identify the vertices of $Q$ with $\mu_v(Q)$, just as we can identify the homology bases of a weave before and after Legendrian mutation.
	
	\begin{remark}
		The crucial difference between algebraic and geometric intersections is captured in the step canceling directed 2-cycles. This cancellation is implemented by default in a quiver mutation, as the arrows of the quiver {\it only} capture algebraic intersections. In contrast, the geometric intersection of homology cycles after a Legendrian mutation will, in general, not coincide with the algebraic intersection. This dissonance will be explored in detail in Section \ref{sec:Proofs}.\hfill$\Box$
	\end{remark}
	
	The following theorem relates the two operations of quiver mutation and Legendrian mutation:
	
	\begin{center}\begin{figure}[h!]{ \includegraphics[width=.8\textwidth]{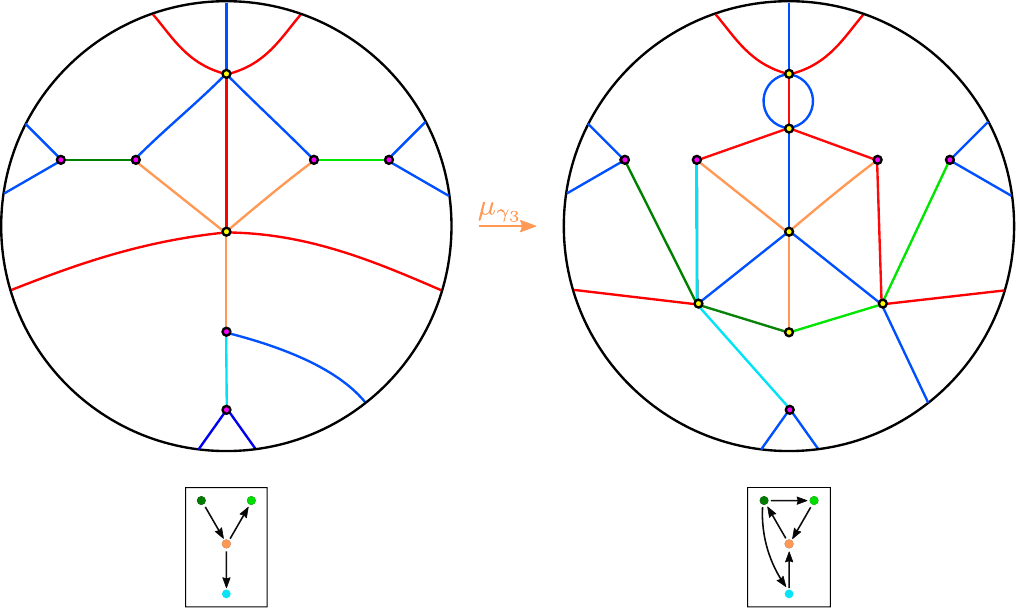}}\caption{Mutation of $\Gamma_0(D_4)$ and its associated intersection quiver at the short {\sf Y}-cycle colored in orange. Note that the sign of the intersection between the dark green {\sf I}-cycle and the orange {\sf Y}-cycle changes from negative to positive, reflecting the reversal of the arrow in the quiver under mutation.}\label{fig:D4Mutation}\end{figure}
	\end{center}
	\begin{theorem}[\cite{CasalsZaslow}, Section 7.3]
		Given a 3-graph $\Gamma$, Legendrian mutation at an embedded cycle $\gamma$ induces a quiver mutation for the associated intersection quivers, taking $Q(\Gamma, \{\gamma_i\})$ to $\mu_\gamma(Q(\Gamma, \{\gamma_i\})).$ \hfill $\Box$
	\end{theorem}
	
	See Figure \ref{fig:D4Mutation} for an example showing the quiver mutation of $Q(\Gamma_0(D_4),\{\gamma_i^{(0)}\})$, $i\in[1,4]$, corresponding to Legendrian mutation applied to $\Lambda(\Gamma_0(D_4)).$ 
	

	\subsection{Microlocal Sheaves and Clusters}

	To introduce the cluster structure mentioned above, we need to define a sheaf-theoretic invariant. We first consider the dg-category of complexes of sheaves of $\C-$modules on $\D^2\times \R$ with constructible cohomology sheaves. For a given 3-graph $\Gamma$ and its associated Legendrian $\Lambda(\Gamma)$, we denote by $\mathcal{C}(\Gamma):=Sh^1_{\Lambda(\Gamma)}(\D^2\times \R)_0$ the subcategory of microlocal rank-one sheaves with microlocal support along $\Lambda(\Gamma)$, which we require to be zero in a neighborhood of $\D^2\times \{-\infty\}$. Here we identify the unit cotangent bundle $T^{\infty, -} ( \D^2\times \R)$ with the first jet space $J^1(\D^2).$ With this identification, the sheaves of $\mathcal{C}(\Gamma)$ are constructible with respect to the stratification given by the Legendrian front $\Sigma(\Gamma).$ Work of Guillermou, Kashiwara, and Schapira implies that that $\mathcal{C}(\Gamma)$ is an invariant under Hamiltonian isotopy 
	\cite{GKS_Quantization}.

\begin{center}\begin{figure}[h!]{ \includegraphics[width=.6\textwidth]{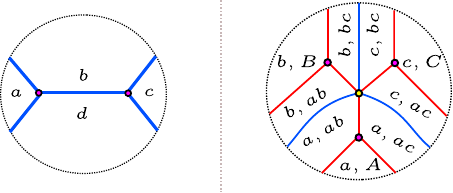}}\caption{The data of the flag moduli space given in the neighborhood of a short {\sf I}-cycle (left) and a short {\sf Y}-cycle (right). Lines are represented by lowercase letters, while planes are written in uppercase. The intersection of the two lines $a$ and $b$ is written as $ab$. } 
			\label{fig:flags}\end{figure}
	\end{center}

	As described in \cite[Section 5.3]{CasalsZaslow}, this category has a combinatorial description. Given a 3-graph $\Gamma$, the data of the moduli space of microlocal rank-one sheaves is equivalent to providing: \begin{enumerate}
		\item[(i)] An assignment to each face $F$ (connected component of $\D^2\backslash G$) of a flag $\mathcal{F}^\bullet(F)$ in the vector space $\C^3$.
		\item[(ii)] For each pair $F_1, F_2$ of adjacent faces sharing an edge labeled by $\sigma_i$, we require that the corresponding flags satisfy 
		$$\mathcal{F}^j(F_1)=\mathcal{F}^j(F_2), \qquad 0\leq j\leq 3, j\neq i, \qquad \text{ and } \qquad \mathcal{F}^i(F_1)\neq \mathcal{F}^i(F_2).$$  
	\end{enumerate}
	Finally, we consider the moduli space of flags satisfying (i) and (ii) modulo the diagonal action of $GL_n$ on $\mathcal{F}^\bullet$. The precise statement \cite[Theorem 5.3]{CasalsZaslow}  is that the flag moduli space, denoted by $\mathcal{M}(\Gamma)$, is isomorphic to the space of microlocal rank-one sheaves $\mathcal{C}(\Gamma)$. 
	Since $\mathcal{C}(\Gamma)$ is an invariant of $\Lambda(\Gamma)$ up to Hamiltonian isotopy, it follows that $\mathcal{M}(\Gamma)$ is an invariant as well. In the {\sf I}-cycle case, when the edges are labeled by $\sigma_1$, the moduli space is determined by four lines $a\neq b \neq c\neq d\neq a$, as pictured in Figure \ref{fig:flags} (left). If the edges are labeled by $\sigma_2$, then the data is given by four planes $A\neq B\neq C \neq D \neq A.$ Around a short {\sf Y}-cycle, the data of the flag moduli space is given by three distinct planes $A\neq B \neq C\neq A $ contained in $\C^3$ and three distinct lines $a\subsetneq A, b\subsetneq B, c\subsetneq C$ with $a\neq b\neq c\neq a,$ as pictured in Figure \ref{fig:flags} (right).
	
\begin{center}\begin{figure}[h!]{ \includegraphics[width=.8\textwidth]{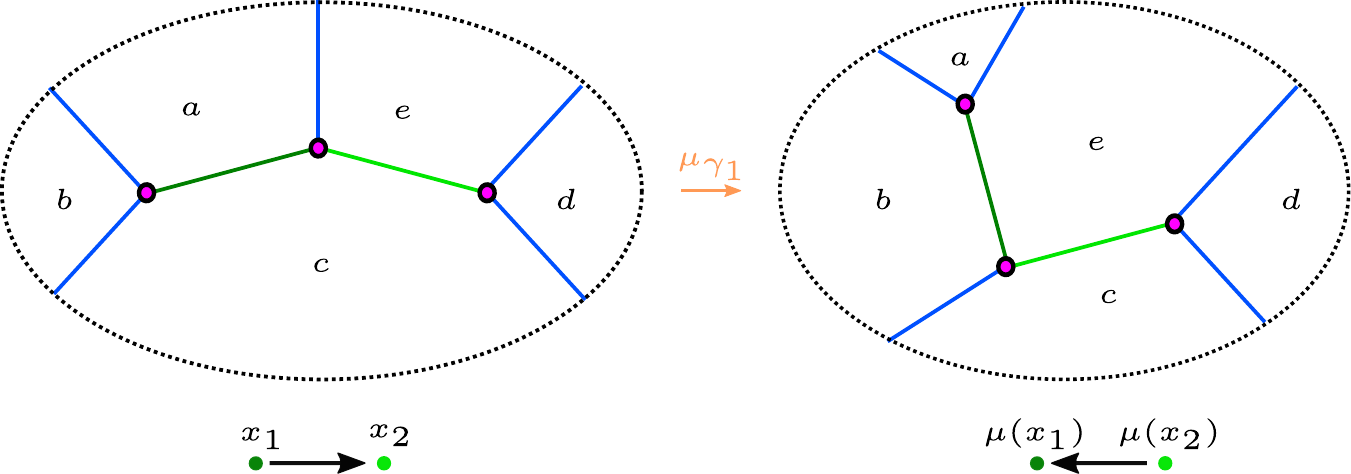}}\caption{Prior to mutating at $\gamma_1,$ we have $\langle\gamma_1,\gamma_2\rangle=-1$. Computing the cross ratios for $\gamma_1$ and $\mu_1(\gamma_1)$ we can see that the cross ratio transforms as $\mu_1(\gamma_1)=\frac{b\wedge c}{c\wedge e}\frac{e\wedge a}{a\wedge b}=x_1^{-1}$ under mutation. Similarly, computing the cross ratios for $\gamma_1$ and $\mu_1(\gamma_2)$ and applying the relation $e\wedge b \cdot a \wedge c= b\wedge c \cdot e\wedge a + a \wedge b \cdot c\wedge e,$ we have $\mu_1(x_2)=\frac{e\wedge a}{a\wedge c}\frac{c\wedge d}{d\wedge e}\left(1 + \frac{a\wedge b}{b\wedge c}\frac{c\wedge e}{e\wedge a}\right).$}\label{fig:clustercalc}\end{figure}\end{center}

	To describe the cluster algebra structure on $\mathcal{C}(\Gamma)$, we need to specify the cluster seed associated to the quiver $Q(\Lambda(\Gamma),\{\gamma_i\})$ via the microlocal monodromy functor $\mu_{mon}$, which is a functor from the category $\mathcal{C}(\Gamma)$ to the category of rank one local systems on $\Lambda(\Gamma)$. As described in \cite{STZ_ConstrSheaves, STWZ}, the functor $\mu_{mon}$ takes a 1-cycle as input and outputs the isomorphism of sheaves given by the monodromy about the cycle. Since it is locally defined, we can compute the microlocal monodromy about an {\sf I}-cycle or {\sf Y}-cycle using the data of the flag moduli space in a neighborhood of the cycle. If we have a short {\sf I}-cycle $\gamma$ with flag moduli space described by the four lines $a, b, c, d$, as in Figure \ref{fig:flags} (left), then the microlocal monodromy about $\gamma$ is given by the cross ratio 
	$$\frac{a\wedge b}{b\wedge c}\frac{c\wedge d}{d\wedge a}$$
	Similarly, for a short {\sf Y}-cycle with flag moduli space given as in Figure \ref{fig:flags} (right), the microlocal monodromy is given by the triple ratio
	$$\frac{B(a)C(b)A(c)}{B(c)C(a)A(b)}$$ 
	As described in \cite[Section 7.2]{CasalsZaslow}, the microlocal monodromy about a 1-cycle gives rise to an $X$-cluster variable at the corresponding vertex in the quiver. Under mutation of the 3-graph, the cross ratio and triple ratio transform as cluster $X$-coordinates. Specifically, if we start with a 3-graph with cluster variables $x_j$, then the cluster variables $x_j'$ of the 3-graph after mutating at $\gamma_i$ are given by the equation 
	$$x_j'=\begin{cases}
		x_j^{-1}  & i=j\\
		x_j(1+x_i^{-1})^{-\langle \gamma_i, \gamma_j\rangle} & \langle \gamma_i, \gamma_j\rangle>0\\
		x_j(1+x_i)^{-\langle \gamma_i, \gamma_j\rangle} & \langle \gamma_i, \gamma_j\rangle<0
	\end{cases}$$ 
	See Figure \ref{fig:clustercalc} for an example.


	The goal of the next section will be to realize each possible mutation of the $D_n$ quiver as a mutation of the corresponding 3-graph. This will imply that there are at least as many exact Lagrangian fillings as cluster seeds of $D_n$-type. There exists a complete classification of all finite mutation type cluster algebras, and in fact, the number of cluster seeds of $D_n$-type is $(3n-2)C_{n-1}$ \cite{FWZ2}.

	\begin{remark} Other than Legendrian weaves, it is not known whether methods of generating exact Lagrangian fillings of $\la(D_n)$ access all possible cluster seeds of $D_n$-type. When constructing fillings of $D_4$ by opening crossings, as in \cite{EHK, Pan-fillings}, experimental evidence suggests that it is only possible to access at most 46 out of the possible 50 cluster seeds by varying the order of the crossings chosen. Of note in the combinatorial setting, we also contrast the 3-graphs $\Gamma(D_4)$ with double wiring diagrams for the torus link $T(3,3)$, which is the smooth type of $\lambda(D_4)$. The moduli of sheaves $\SC(\Gamma(D_4))$ for $\Gamma(D_4)$ embeds as an open positroid cell into the Grassmanian $Gr(3,6)$ \cite{CasalsGao}, so we can identify some cluster charts with double wiring diagrams. The double wiring diagrams associated to $Gr(3,6)$ only access 34 out of 50 distinct cluster seeds via local moves applied to an initial double wiring diagram \cite{FWZ1}. \hfill $\Box$ 
	\end{remark}


	\section{Proof of Main Results}\label{sec:Proofs}

	In this section, we state and prove Theorem \ref{technical}, which implies Theorem \ref{main}. 
	The following definitions relate the algebraic intersections of cycles to geometric intersections in the context of 3-graphs.

	\begin{definition}
		A 3-graph $\Gamma$ with associated basis $\{\gamma_i\},$ $i\in [1, b_1(\Lambda(\Gamma)]$ of $H_1(\Lambda(\Gamma); \Z)$ is \emph{sharp at a cycle} $\gamma_j$ if, for any other cycle $\gamma_k\in \{\gamma_i\}$, the geometric intersection number of $\gamma_j$ with $\gamma_k$ is equal to the algebraic intersection $\langle \gamma_j, \gamma_k\rangle$. 
		
		$\Gamma$ is \emph{locally sharp} if, for any cycle $\gamma\in \{\gamma_i\},$ there exists a sequence of Legendrian Surface Reidemeister moves 
		taking $\Gamma$ to some other 3-graph $\Gamma'$ such that $\Gamma'$ is sharp at the corresponding cycle $\gamma'\in H_1(\Lambda(\Gamma');\Z)$.
		
		A 3-graph $\Gamma$ with a set of cycles $\Gamma$ is  \emph{sharp} if $\Gamma$ is sharp at all $\gamma_i\in \{\gamma_i\}$. \hfill $\Box$
	\end{definition}
	
	
	For 3-graphs that are not sharp, it is possible that a sequence of mutations will cause a cycle to become immersed. This is the only obstruction to weave realizability. Therefore, sharpness is a desirable property for our 3-graphs, as it simplifies our computations and helps us avoid creating immersed cycles. We will not be able to ensure sharpness for all $\Gamma(D_n)$ that arise as part of our computations, (e.g., see the type III.i normal form in Figure \ref{fig:normalforms}) but we will be able to ensure that each of our 3-graphs is locally sharp.
	
	
	\subsection{Proof of Theorem \ref{main}}

	The following result is slightly stronger than the statement of Theorem \ref{main}, as we are able to show that each 3-graph in our sequence of mutations is locally sharp.
	
	
	\begin{theorem}\label{technical}
		Let $\mu_{v_1}, \dots, \mu_{v_k}$ be a sequence of quiver mutations, with initial quiver $Q(\Gamma_0(D_n), \{\gamma_i^{(0)}\})$. Then, there exists a sequence $\Gamma_0(D_n), \dots , \Gamma_k(D_n)$ of 3-graphs such that 
		\begin{enumerate}
			\item[i.] $\Gamma_{j-1}(D_n)$ is related to $\Gamma_j(D_n)$ by mutation at a cycle $\gamma_j$ and by Legendrian Surface Reidemeister moves I, II, and III. The cycle $\gamma_j$ represents the vertex $v_j$ in the intersection quiver and it is given by one of the cycles in the initial basis $\{\gamma_i^{(0)}\}$ after mutation and Reidemeister moves.\\
			
			\item[ii.] $\Gamma_j(D_n)$ is sharp at $\gamma_j$.\\
			
			\item[iii.] $\Gamma_j(D_n)$ is locally sharp.\\
			
			\item[iv.] The basis of cycles for $\Gamma_j(D_n)$, obtained from the initial basis $\{\gamma_i^{(0)}\}$ by mutation and Reidemeister moves, consists entirely of short {\sf Y}-cycles and short {\sf I}-cycles.
		\end{enumerate}
	\end{theorem}
	
	The conditions ii-iv allow us to continue to iterate mutations after applying a small number of simplifications at each step. Theorem \ref{main} thus follows from Theorem \ref{technical}.
	
	\bigskip
	
	\begin{proof}

		We proceed by organizing the 3-graphs arising from any sequence of mutations of $\Gamma_0(D_n)$ into four types, in line with the organization scheme introduced by Vatne for quivers of $D_n$-type \cite{Vatne}. Vatne's classification of quivers in the mutation class of $D_n$-type uses the configuration of a certain subquiver to define the different types. Outside of that subquiver, there are a number of disjoint subquivers of $A_n$-type that are referred to as $A_n$ tail subquivers. We will refer to the corresponding cycles in the 3-graph as $A_n$ tail subgraphs, or simply $A_n$ tails when it is clear from context whether we are referring to the quiver or the 3-graph. For each type, Vatne describes the results of quiver mutation at different vertices, which can depend on the existence of $A_n$ tail subquivers. See Figures \ref{fig:1Quivers}, \ref{fig:IIQuivers}, \ref{fig:3Quivers}, and \ref{fig:4Quivers} for the four types and their mutations.

		\textbf{Notation.} As mentioned in the previous section, cycles are pictured as colored edges for the sake of visual clarity. Throughout this section, we denote all of the dark green cycles by $\gamma_1,$ light green cycles by $\gamma_2$, orange cycles by $\gamma_3$, light blue cycles by $\gamma_4$, pink cycles by $\gamma_5$, purple cycles by $\gamma_6$, and olive cycles by $\gamma_7$. With this notation, $\gamma_i$ will correspond to the vertex labeled by $v_i$ in the quivers given below. 
		
		\textbf{$\boldsymbol{A_n}$ Tails.} 
		We briefly describe the behavior of the $A_n$ tail subquivers, as given in \cite{Vatne}, in terms of weaves. Any of the $n$ vertices in an $A_n$ tail subquiver can have valence between 0 and 4. Cycles in the quiver are oriented with length 3. 
		If a vertex $v$ has valence 3, then two of the edges form part of a 3-cycle, while the third edge is not part of any 3-cycle. If $v$ has valence 4, then two of the edges belong to one 3-cycle and the remaining two edges belong to a separate 3-cycle.

		Any $A_n$ tail of the quiver can be represented by a sharp configuration of $n$ {\sf I}-cycles in the 3-graph. See Figure \ref{fig:Tails} for an identification of {\sf I}-cycles with quiver vertices of a given valence. Mutation at any vertex $v_i$ in the quiver corresponds to mutation at the {\sf I}-cycle $\gamma_i$ in the 3-graph, so it is readily verified that mutation preserves the number of {\sf I}-cycles and requires no application of Legendrian Surface Reidemeister moves to simplify. The sequences of mutations given in the remainder of the proof   As a consequence, any sequence of $A_n$ tail mutations is weave realizable, and a sharp 3-graph remains sharp after mutation at $A_n$ tail {\sf I}-cycles that only intersect other $A_n$ tail {\sf I}-cycles. 
		\begin{center}\begin{figure}[h!]{ \includegraphics[width=.8\textwidth]{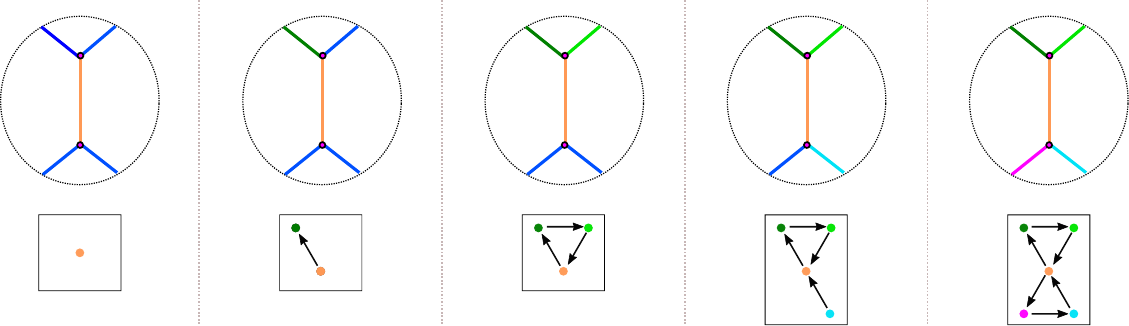}}\caption{All possible arrangements of {\sf I}-cycles in an $A_n$ tail of the 3-graph corresponding to a given vertex in the $A_n$ tail subquiver of valence between 0 and 4.}\label{fig:Tails}\end{figure}\end{center}

		
		\textbf{Normal Forms.} 
		For each of the four types of $D_n$ quivers described in \cite{Vatne}, we give a set of specific subgraphs of $\Gamma(D_n)$, which we refer to as normal forms. These normal forms are pictured in Figure \ref{fig:normalforms}. 
		We indicate the possible existence of $A_n$ tail subgraphs by an unfilled circle. In our discussion below, we will say that an edge of the 3-graph carries a cycle if it is part of a homology cycle. We will generally use this terminology to specify which edges cannot carry a cycle.  
		
		
		

		\begin{center}\begin{figure}[h!]{ \includegraphics[width=\textwidth]{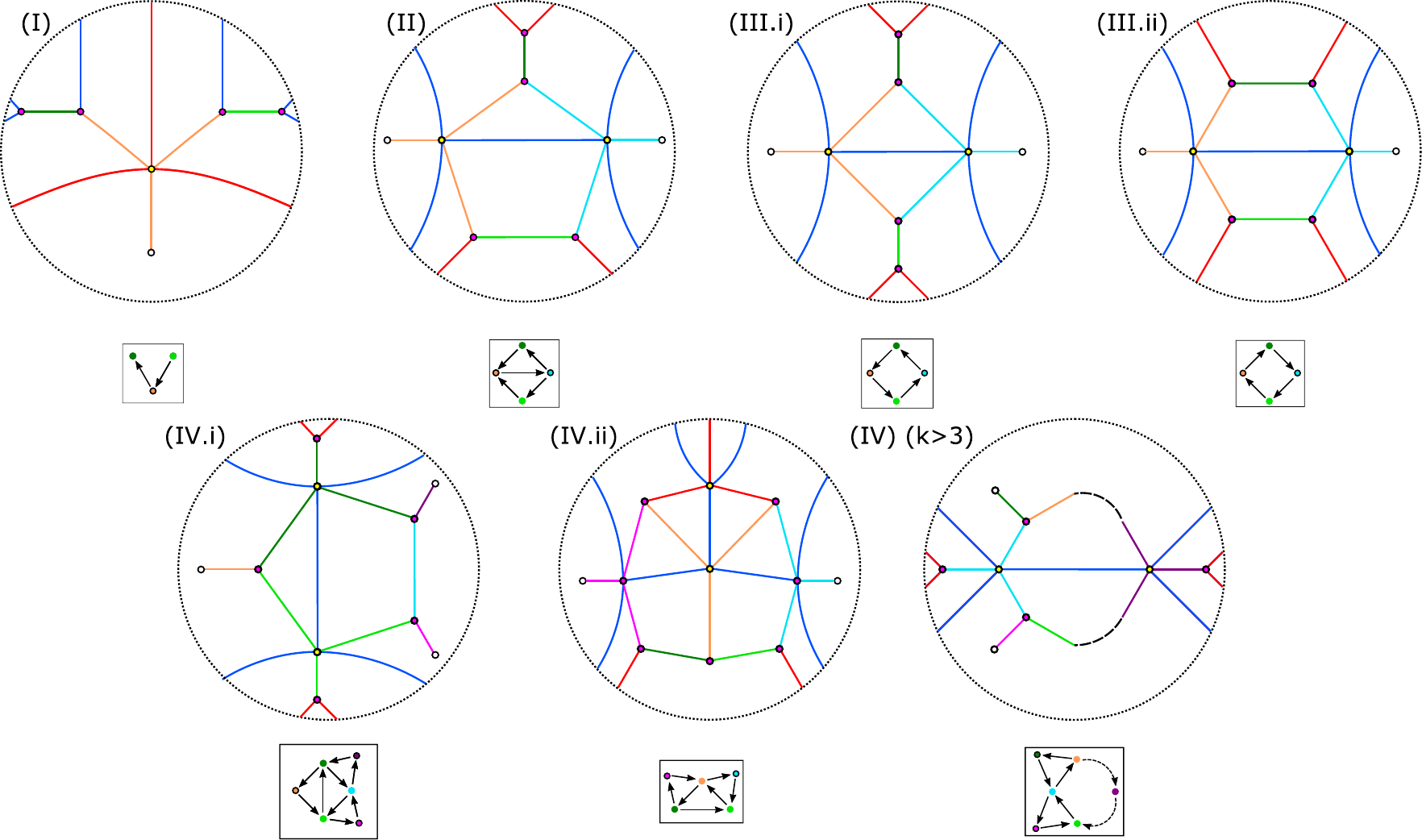}}\caption{Normal forms labeled by their type. The possible addition of {\sf I}-cycles corresponding to $A_n$ tails of the quiver are represented by unfilled circles appended to the end of edges that do not intersect the boundary. }\label{fig:normalforms}\end{figure}\end{center}

		For each possible quiver mutation, we describe the possible mutations of the 3-graph and show that the result matches the quiver type and retains the properties listed in Theorem \ref{technical} above. In addition, the Legendrian Surface Reidemeister moves we describe ensure that the $A_n$ tail subgraphs continue to consist solely of short {\sf I}-cycles. If the mutation results in a long {\sf I}-cycle or pair of long {\sf I}-cycles connecting our $A_n$ tail to the rest of the 3-graph, we can simplify by applying a sequence of $n$ push-throughs to ensure that these are all short {\sf I}-cycles. It is readily verified that we can always do this and that no other simplifications of the $A_n$ tails are required following any other mutations. We include $A_n$ tail cycles only where relevant to the specific mutation. In our computations below, we generally omit the final steps of applying a series of push-throughs to make any long {\sf I} or {\sf Y}-cycles into short {\sf I} or {\sf Y}-cycles. Figure \ref{fig:push-throughs} provides an example where these push-throughs are shown for both an {\sf I}-cycle and a {\sf Y}-cycle. 

		In order to simplify the overall presentation of the normal forms and the computations below, we allow for the following variations in the Type I and Type IV cases. In the Type I case,
		mutating at either of the short {\sf I}-cycles $\gamma_1$ or $\gamma_2$ in the Type I normal form produces one of four possible configurations of the cycles $\gamma_1, \gamma_2,$ and $\gamma_3$ in a  3-graph corresponding to a Type I quiver. Since these mutations are readily computed, we simplify our presentation by giving a single normal form rather than four, and describing the relevant mutations of two of the four possible 3-graphs in figures \ref{fig:1to1}, \ref{fig:1to2}, \ref{fig:1to4.1}, and \ref{fig:1to4.2}. The remaining cases can be seen by swapping the cycle(s) to the left of the short {\sf Y}-cycle with the cycle(s) to the right of it. This symmetry corresponds to reversing all of the arrows in the quiver. In general, we will implicitly appeal to similar symmetries of the normal form 3-graphs to reduce the number of cases we must consider. In the Type IV case, the edge(s) corresponding to $\gamma_3, \gamma_5$ or $\gamma_6$ need not carry a cycle. See the discussion of Type IV quiver mutations below for a more detailed description.
		

		

		\textbf{Type I.} 
		We start with 3-graphs, always endowed with a homology basis, whose associated intersection quivers are a Type I quiver. See Figure \ref{fig:1Quivers} for the relevant quiver mutations. 
		
		\begin{center}\begin{figure}[h!]
				{\includegraphics[width=\textwidth]{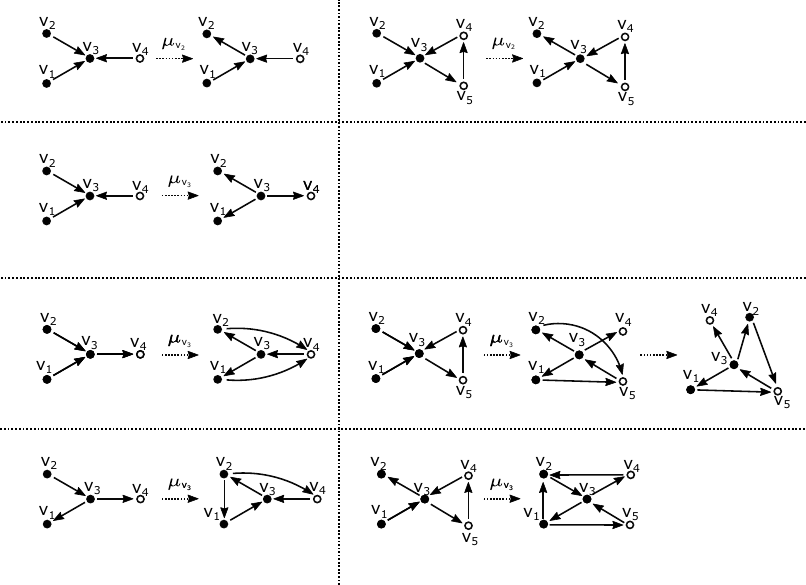}}
				\caption{From top to bottom, two Type I to Type I quiver mutations, Type I to Type II quiver mutations, and Type I to Type IV quiver mutations. The arrow labeled by $\mu_{v_i}$ indicates mutation at the vertex $v_i$. Unfilled circles represent potential $A_n$ tails. In each line, the first quiver mutation shows the case where $v_3$ is only adjacent to one $A_n$ tail vertex, while the second quiver mutation shows the case where $v_3$ is adjacent to two $A_n$ tail vertices. Note that reversing the direction of all of the arrows simultaneously before mutating gives additional possible quiver mutations of the same type.
				}\label{fig:1Quivers}
			\end{figure}
		\end{center}

		\begin{center}\begin{figure}[h!]
				{\includegraphics[width=\textwidth]{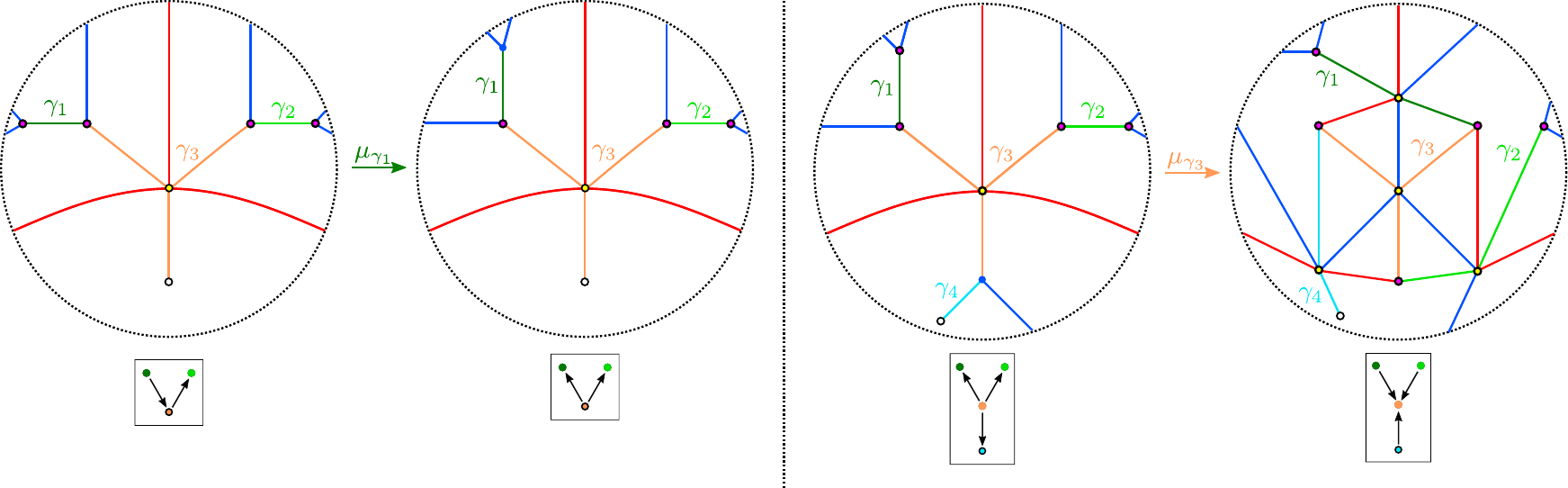}}
				\caption{Type I to Type I mutation. Arrows labeled by $\mu$ indicate mutation at a cycle of the same color.}\label{fig:1to1}
			\end{figure}
		\end{center}
		
		\begin{itemize}
			\item[i.] (Type I to Type I) There are two possible Type I to Type I mutations of 3-graphs depicted in Figure \ref{fig:1to1} (left) and (right). As shown in Figure \ref{fig:1to1} (left), mutation at $\gamma_1$ only affects the sign of the intersection of $\gamma_1$ with the $\gamma_3$. This reflects the fact that the corresponding quiver mutation has only reversed the orientation of the edge between $v_1$ and $v_3$. Mutating at any other {\sf I}-cycle is equally straightforward and yields a Type I to Type I mutation as well.\\

			\item[ii.] (Type I to Type I) For the second possible Type I to Type I mutation, we proceed as pictured in Figure \ref{fig:1to1} (right).  Mutation at $\gamma_3$ does not create any new additional geometric or algebraic intersections. Instead, it takes positive intersections to negative intersections and vice versa. This is reflected in the quivers pictured underneath the 3-graphs, as the orientation of edges has reversed under the mutation. As explained above, we could simplify the resulting 3-graph by applying a push-through move to each of the long {\sf I}-cycles to get a sharp 3-graph where the homology cycles are made up of a single short {\sf Y}-cycle and some number of short {\sf I}-cycles.\\ 
			
			\item[iii.] (Type I to Type II) In Figure \ref{fig:1to2} we consider the cases where the {\sf Y}-cycle $\gamma_3$ intersects one {\sf I}-cycle (top) or two {\sf I}-cycles (bottom) in the $A_n$ tail subgraph. Mutation at $\gamma_3$ introduces an intersection between $\gamma_2$ and $\gamma_4$ that causes the second 3-graph in of each mutation sequences to no longer be sharp. Applying a push-through 
			to $\gamma_2$ resolves this intersection so that the geometric intersection between $\gamma_2$ and $\gamma_4$ matches their algebraic intersection. This simplification ensures that the result of $\mu_{\gamma_3}$ is a sharp 3-graph that matches the Type II normal form. If we compare the mutations in the top and bottom sequences, we can see that the presence of the $A_n$ tail cycle $\gamma_5$ does not affect the computation.\\
			
			\begin{center}
				\begin{figure}[h!]
					{\includegraphics[width=.9\textwidth]{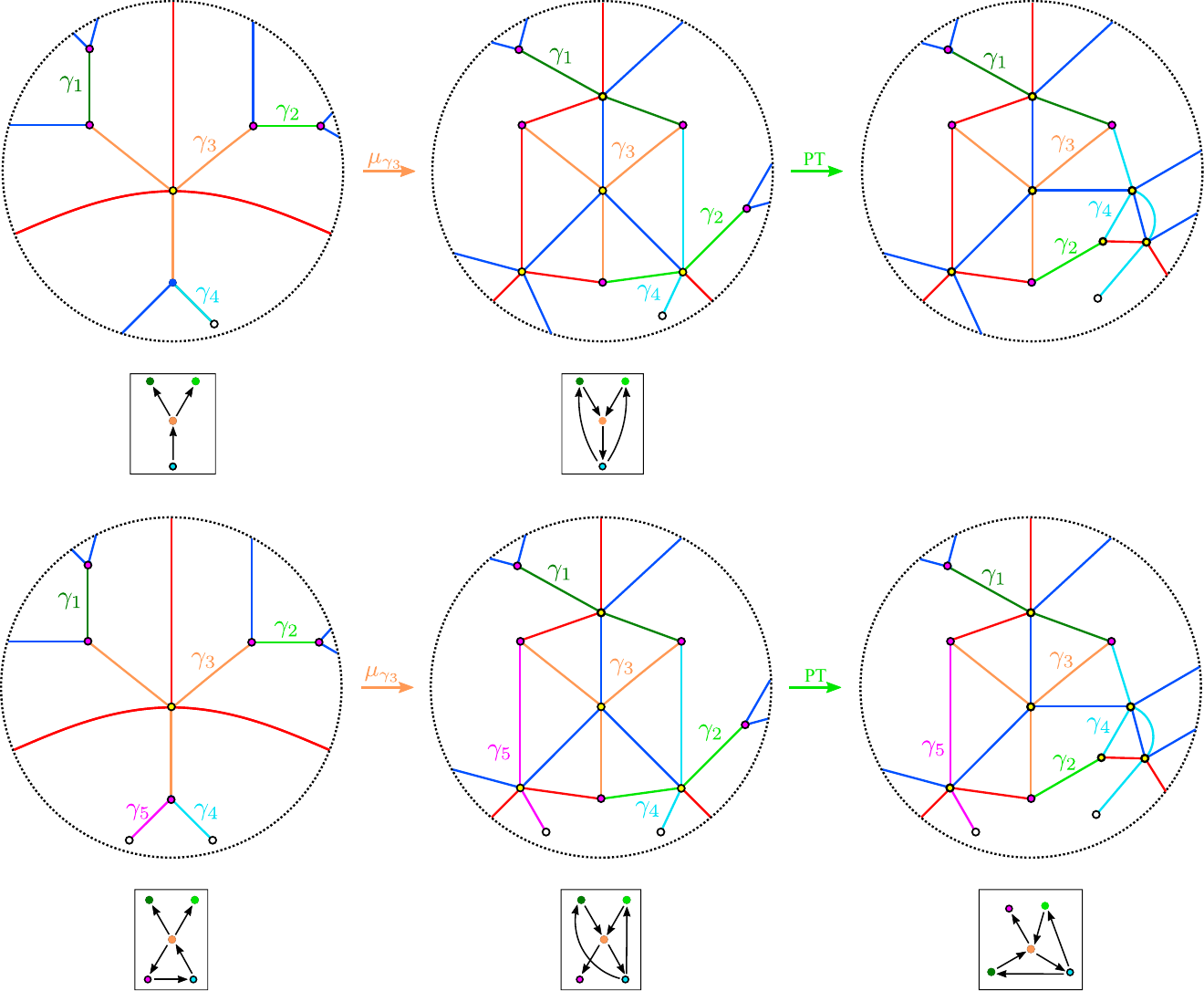}}\caption{Type I to Type II mutations. Legendrian Surface Reidemeister are moves labeled as in Theorem 2, Figure \ref{fig:Moves}.}\label{fig:1to2}
				\end{figure}
			\end{center}
			\item[iv.] (Type I to Type IV.i) We now consider the first of two Type I to Type IV mutations, shown in Figure \ref{fig:1to4.1}. Starting with the configuration of cycles at the left of each sequence and mutating at $\gamma_3$ causes $\gamma_1$ and $\gamma_2$ to cross. Applying a push-through to $\gamma_1$ or to $\gamma_2$ (not pictured) simplifies the resulting intersection and yields a Type IV.i normal form made up of the cycles $\gamma_1, \gamma_2, \gamma_3, $ and $\gamma_4$. The sequences on the top and bottom of Figure \ref{fig:1to4.1} differ only by the presence of the $A_n$ tail cycle $\gamma_5.$\\
			
			\begin{center}
				\begin{figure}[h!]
					{\includegraphics[width=.9\textwidth]{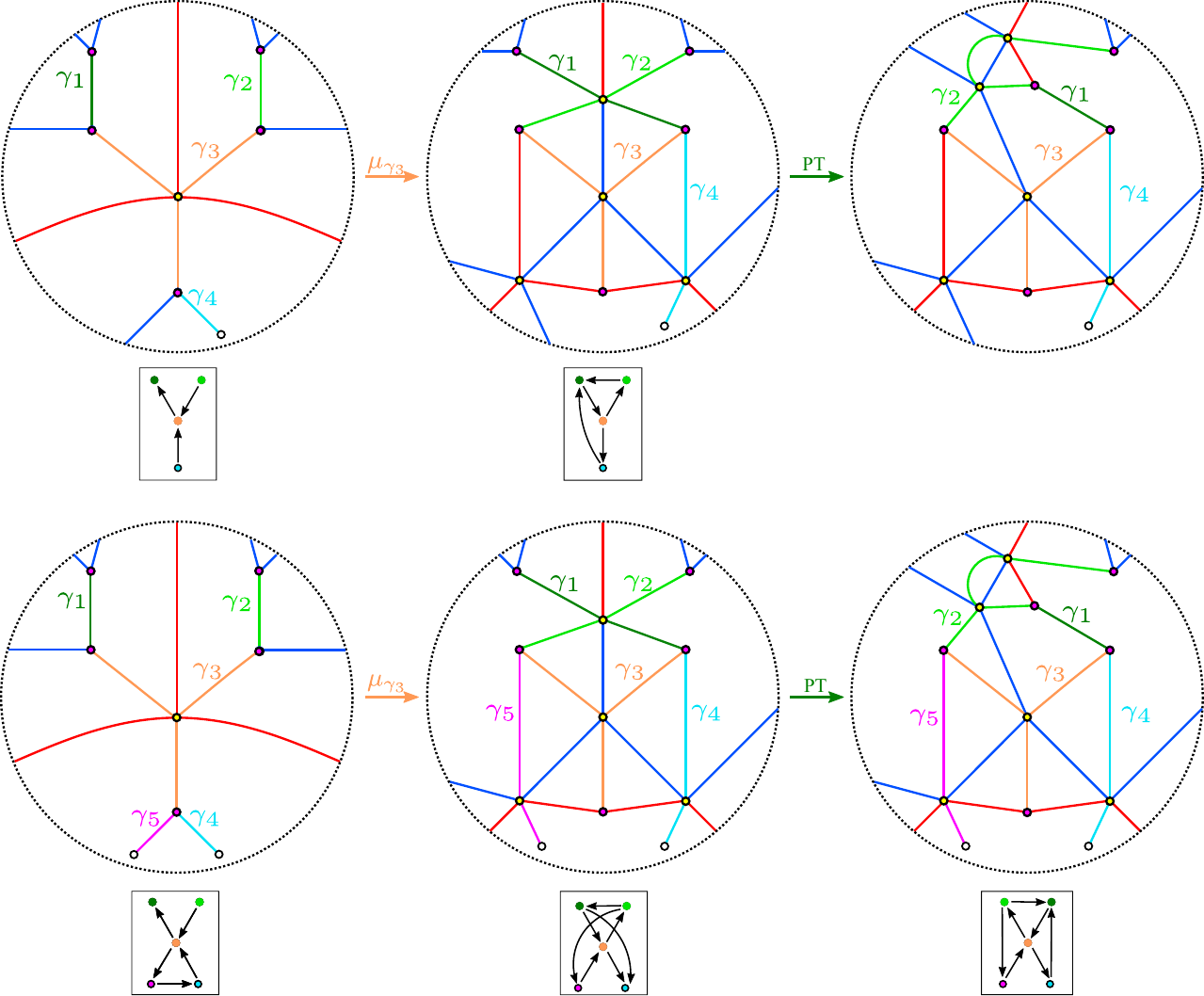}}\caption{Type I to Type IV.i mutations.}\label{fig:1to4.1}
				\end{figure}
			\end{center}
			
			\item[v.] (Type I to Type IV.ii) In Figure \ref{fig:1to4.2},  we consider the cases where $\gamma_1$ intersects one {\sf I}-cycle (top) or two {\sf I}-cycles (bottom) in the $A_n$ tail subgraph, as we did in the Type I to Type II case. 
			As in the Type I to Type II case, we must apply a push-through to resolve the new intersections between that cause the second 3-graph in each sequence to fail to be sharp. When we include both $\gamma_4$ and $\gamma_5$ in the sequence on the right, we get two new intersections after mutating, and therefore require two push-throughs. Note that in the IV.ii case, we must first apply the push-through to $\gamma_1$ and $\gamma_2$ in order to ensure that we can apply a push-through to any additional cycles in the $A_n$ tail subgraph. This causes the {\sf Y}-cycles of the graph to correspond to different vertices in the quiver than in the Type IV.i normal form, which is the main reason we distinguish between the normal forms for Type IV.i and Type IV.ii. \\
			
		\end{itemize}

		\begin{center}\begin{figure}[h!]
				{\includegraphics[width=\textwidth]{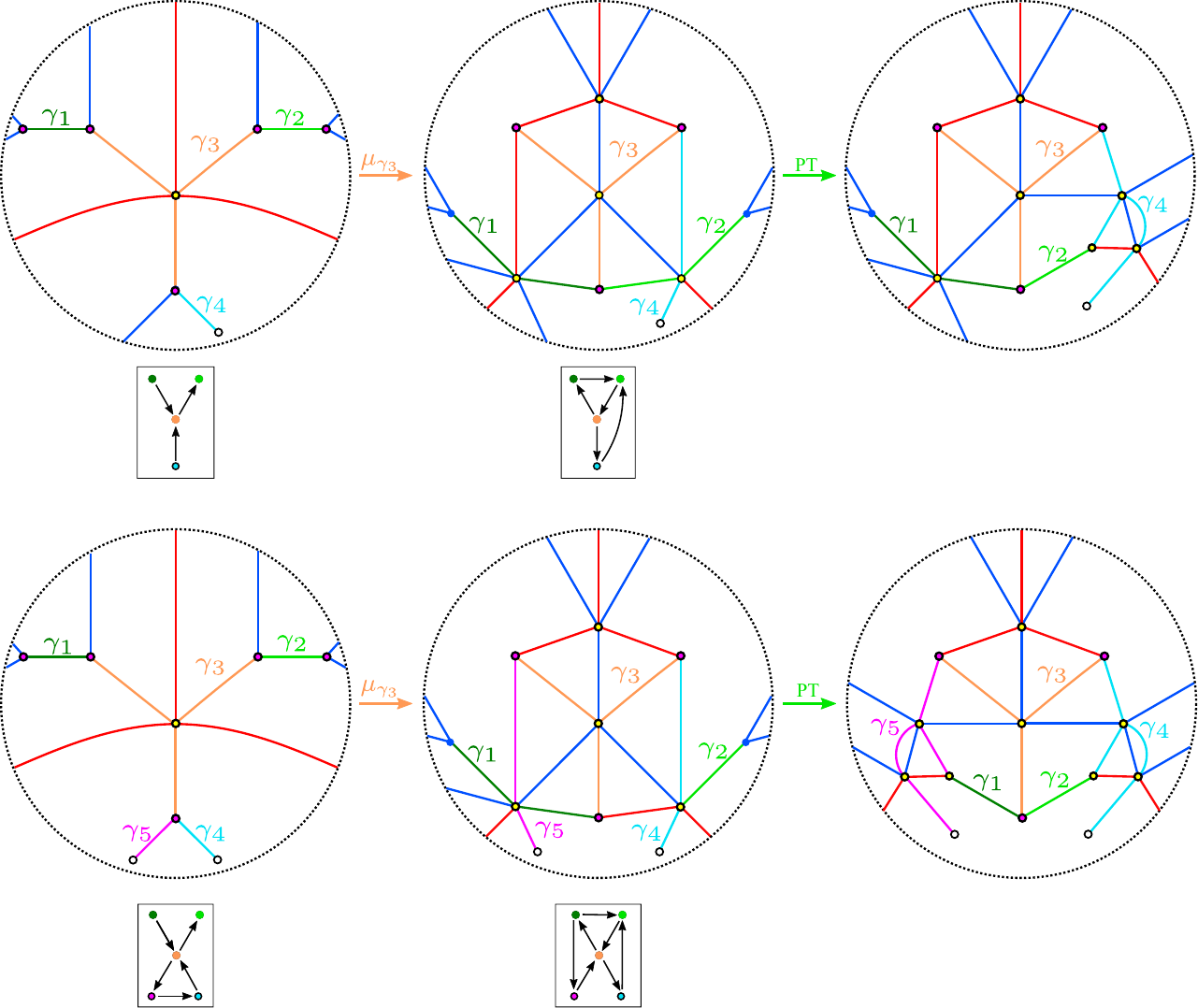}}\caption{Type I to Type IV.ii mutations. }\label{fig:1to4.2}
			\end{figure}
		\end{center}
		In Figure \ref{fig:push-throughs} we show how to apply push-throughs to completely simplify the long {\sf I}- and {\sf Y}-cycles pictured in the Type I to Type IV.ii graph. As mentioned above, these push-throughs are identical to any other computation required to simplify our resulting 3-graphs to a set of short
		{\sf I}- and {\sf Y}-cycles. 
		
		\begin{center}\begin{figure}[h!]
				{\includegraphics[width=\textwidth]{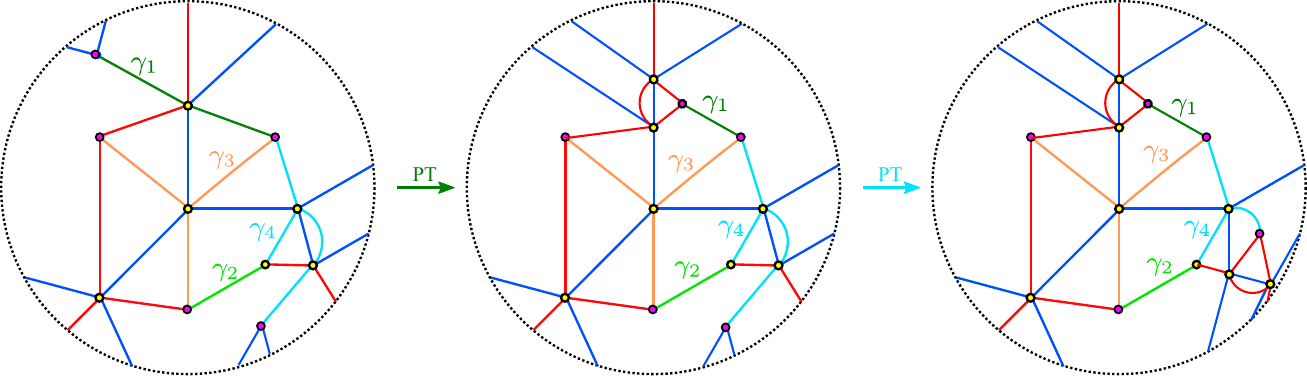}}\caption{Push-through examples. The first push-through move simplifies the long {\sf I}-cycle labeled $\gamma_1$, while the second simplifies the long {\sf Y}-cycle labeled $\gamma_4$.}\label{fig:push-throughs}
			\end{figure}
		\end{center}
		The above cases describe all possible mutations of the Type 1 normal form. Each of these mutations yields a sharp 3-graph with short {\sf I}-cycles and {\sf Y}-cycles, as desired.

		\textbf{Type II.} We now consider mutations of our Type II normal form. See Figure \ref{fig:IIQuivers} for the relevant quivers.	As shown in the figure, performing a quiver mutation at the 2-valent 
		vertices labeled by $v_1$ or $v_2$ yields a Type III quiver, while a quiver mutation at the vertices labeled $v_3$ or $v_4$ yields either another Type II quiver or a Type I quiver, depending on the intersection of $v_3$ or $v_4$ with any $A_n$ tail subquivers.  
		\begin{center}
			\begin{figure}[h!]
				{\includegraphics[width=\textwidth]{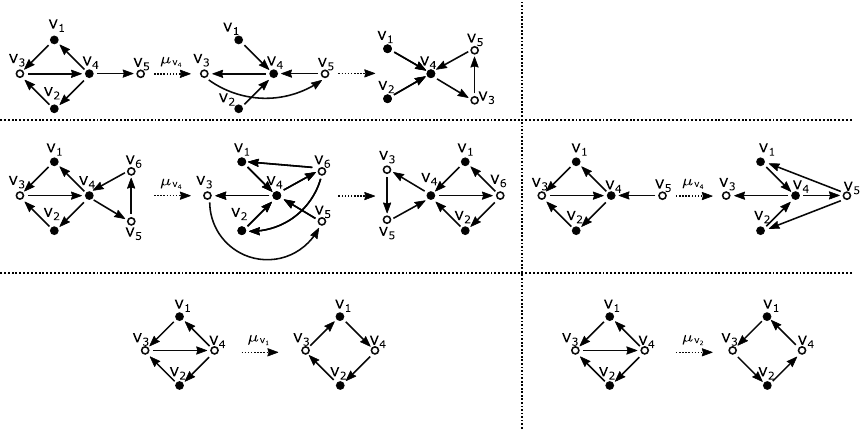}}\caption{From top to bottom, Type II to Type I, Type II to Type II, and Type II to Type III quiver mutations.}\label{fig:IIQuivers}
			\end{figure}
		\end{center}

		\begin{center}
			\begin{figure}[h!]
				{\includegraphics[width=.9\textwidth]{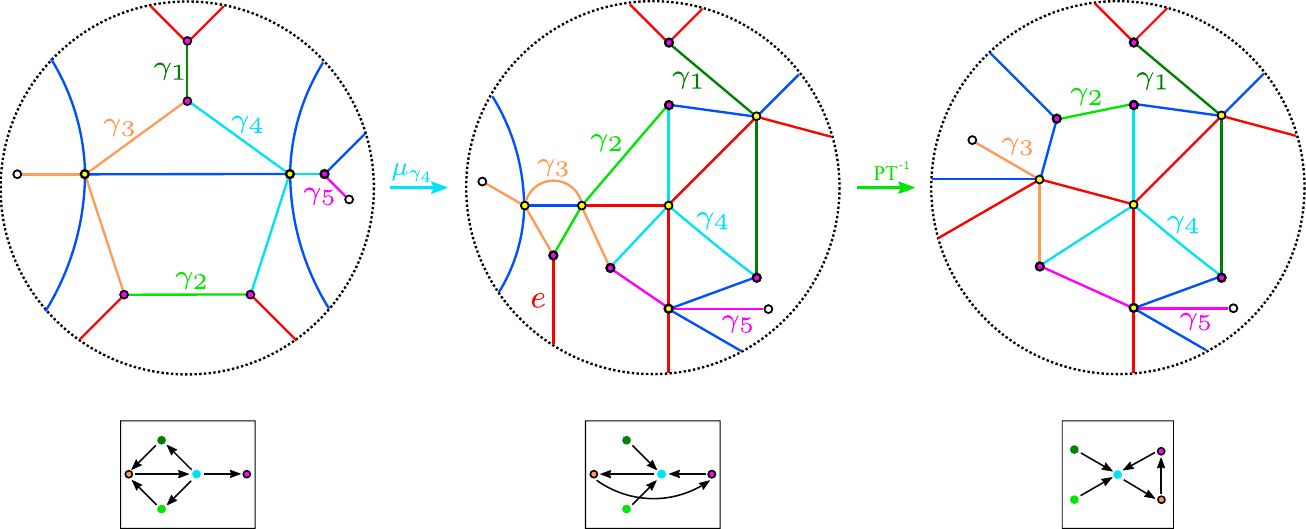}}\caption{Type II to Type I mutations. The red $e$ labels an edge in the 3-graph that does not carry a cycle.
				}\label{fig:2to1}
			\end{figure}
		\end{center}

		
		\begin{itemize}
			\item[i.] (Type II to Type I) We first consider the sequence of 3-graphs pictured in Figure \ref{fig:2to1}. Mutation at $\gamma_4$ results in a new geometric intersection between $\gamma_2$ and $\gamma_3$ even though their algebraic intersection is zero. We can resolve this by applying a reverse push-through at the trivalent vertex where $\gamma_2$ and $\gamma_3$ meet. The resulting 3-graph is sharp, as $\gamma_2$ and $\gamma_3$ no longer have any geometric intersection. This computation is identical if $\gamma_3$ were to intersect a single $A_n$ tail cycle and we mutated at $\gamma_3$ instead. Note that here we require the red edge adjacent labeled $e$ to not carry a cycle, as specified by our normal form.\\
			
			\begin{center}
				\begin{figure}[h!]
					{\includegraphics[width=\textwidth]{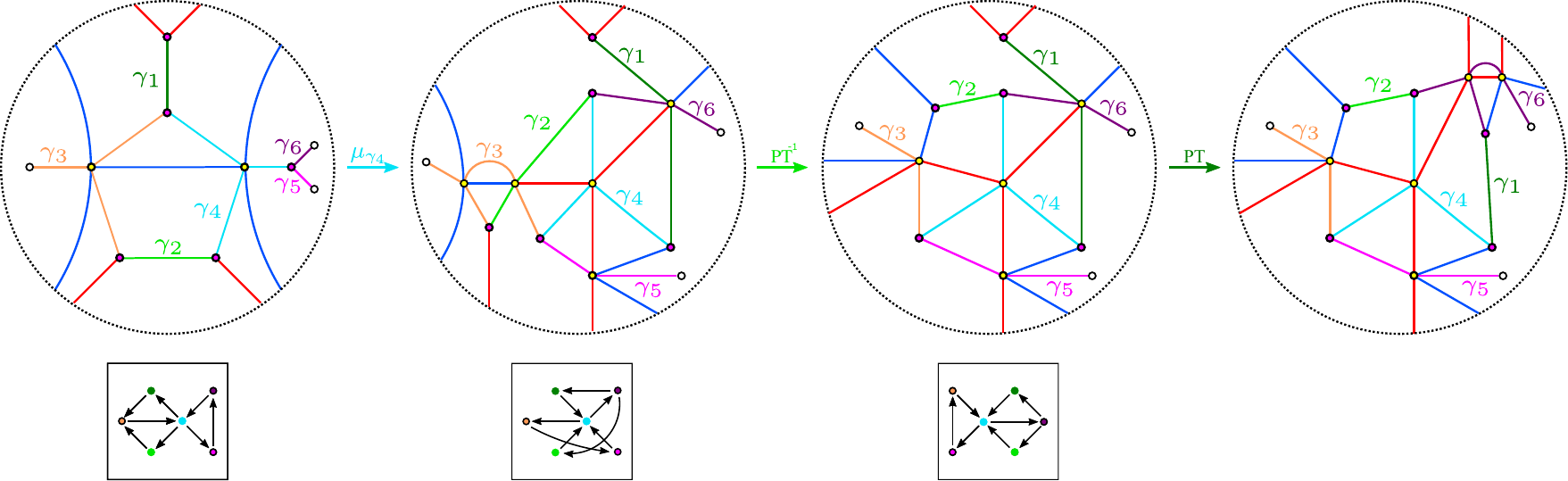}}\caption{Type II to Type II mutations.
					}\label{fig:2to2}
				\end{figure}
			\end{center}

			\item[ii.] (Type II to Type II) We now consider the sequence shown in Figure \ref{fig:2to2}. After mutating at $\gamma_4$, we have the same intersection between $\gamma_2$ and $\gamma_3$ as in the previous case. We again resolve this intersection by applying a reverse push-through at the same trivalent vertex. In this case, we also have an intersection between $\gamma_1$ and $\gamma_6,$ which we resolve via push-through of $\gamma_1$. As a result, $\gamma_6$ becomes a {\sf Y}-cycle, and the Type II normal form is now made up of the cycles $\gamma_1,$ $\gamma_2$, $\gamma_4,$ and $\gamma_6$, while $\gamma_3$ becomes an $A_n$ tail cycle. \\
			

			\begin{center}
				\begin{figure}[h!]
					{\includegraphics[width=\textwidth]{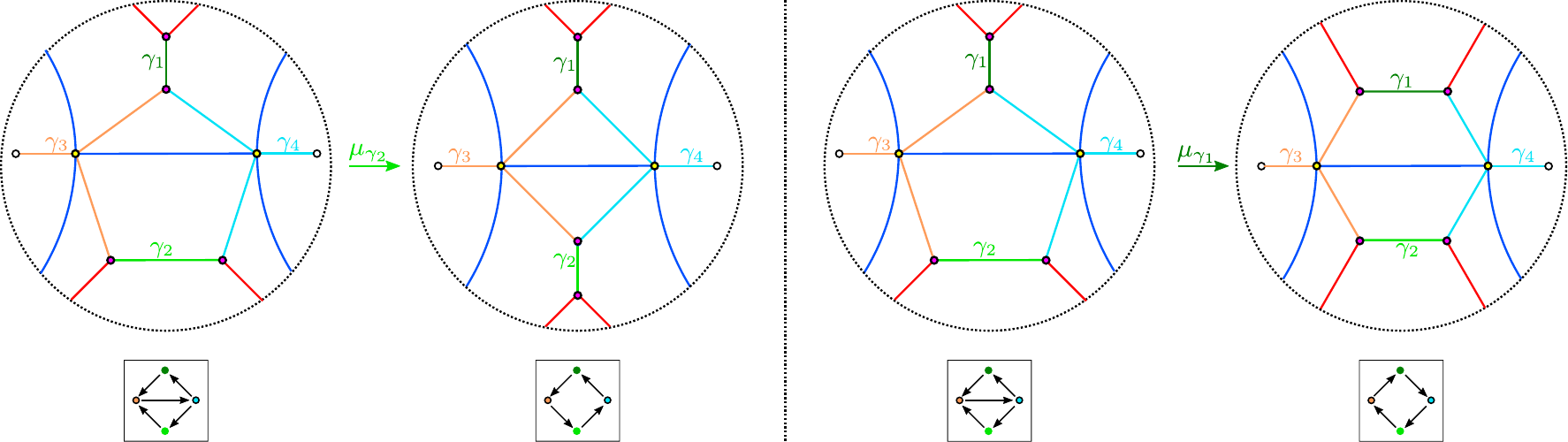}}\caption{Type II to Type III mutations.
					}\label{fig:2to3}
				\end{figure}
			\end{center}
			
			\item[iii.] (Type II to Type III.i) Mutation at $\gamma_1$ or $\gamma_2$ in the Type II normal form yields either of the Type III normal forms. 
			In the sequence on the left of Figure \ref{fig:2to3}, mutation at $\gamma_2$ leads to a geometric intersection between $\gamma_3$ and $\gamma_4$ at two trivalent vertices. Since the signs of these two intersections differ, the algebraic intersection $\langle \gamma_3, \gamma_4\rangle$ is zero, so the resulting 3-graph is not sharp. However, it is sharp at $\gamma_1$ and $\gamma_2$, and applying a flop to the 3-graph removes the geometric intersection between $\gamma_3$ and $\gamma_4$ at the cost of introducing the same intersection between $\gamma_1$ and $\gamma_2$. Therefore, applying the flop does not make the 3-graph sharp, but it does show that the 3-graph resulting from our mutation is locally sharp at every cycle.\\
			
			\item[iv.] (Type II to Type III.ii) In the sequence on the right of Figure \ref{fig:2to3}, mutation at $\gamma_1$ yields a sharp 3-graph that matches the Type III.ii normal form.\\
		\end{itemize}

		\textbf{Type III:} Figure \ref{fig:3Quivers} illustrates the Type III quiver mutations. Figures \ref{fig:3to2}, \ref{fig:3.1to4}, and \ref{fig:3.2to4} depict the corresponding Legendrian mutations of the Type III normal forms.

		\begin{center}
			\begin{figure}[h!]
				{\includegraphics[width=.8\textwidth]{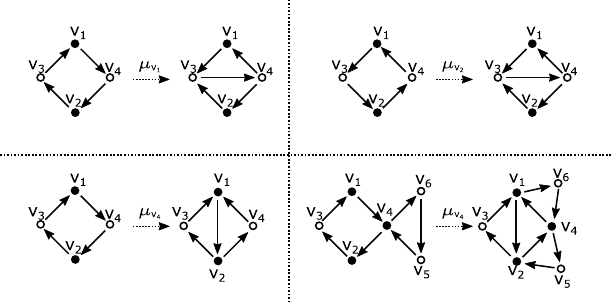}}\caption{Type III to Type II quiver mutations (top) and Type III to Type IV quiver mutations (bottom).}\label{fig:3Quivers}
			\end{figure}
		\end{center}

		\begin{center} 
			\begin{figure}[h!]
				{\includegraphics[width=\textwidth]{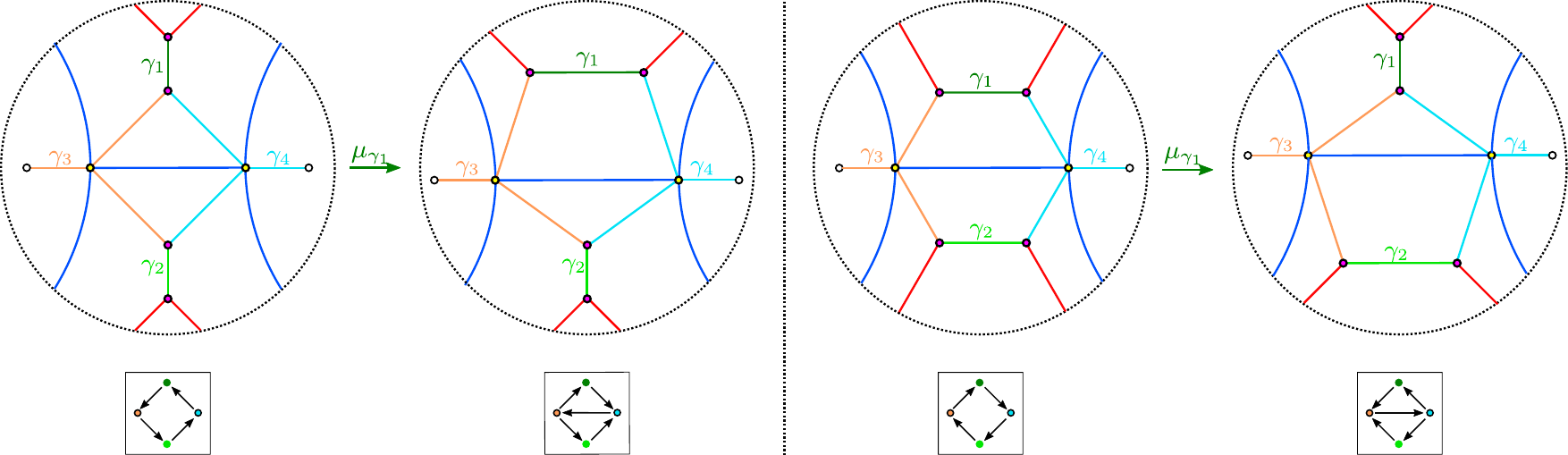}}\caption{Type III.i to Type II mutations (left) and Type III.ii to Type II mutations (right).}\label{fig:3to2}
			\end{figure}
		\end{center}

		\begin{itemize}
			\item[i.] (Type III.i to Type II) We first consider the sequence of 3-graphs in Figure \ref{fig:3to2} (left). Mutating at $\gamma_1$ or $\gamma_2$ immediately yields a Type II normal form. Mutating at $\gamma_1$ and $\gamma_2$ in succession yields a Type III.ii normal form. Note that if the 3-graph were not sharp at $\gamma_1$ or $\gamma_2$ we would first need to apply a flop. We can always apply this move because the 3-graph is locally sharp at each of its cycles. See the Type III.i to Type IV.i subcase below for an example where we demonstrate this move.\\

			\item[ii.] (Type III.ii to Type II) In the sequence on the right of Figure \ref{fig:3to2}, mutation at either $\gamma_1$ or $\gamma_2$ yields a Type II normal form. Mutation at $\gamma_1$ and $\gamma_2$ in succession yields a Type III.i normal form. Therefore, applying these two moves in succession can take us between both of our Type III normal forms.\\
			
			\begin{center} 
				\begin{figure}[h!]
					{\includegraphics[width=\textwidth]{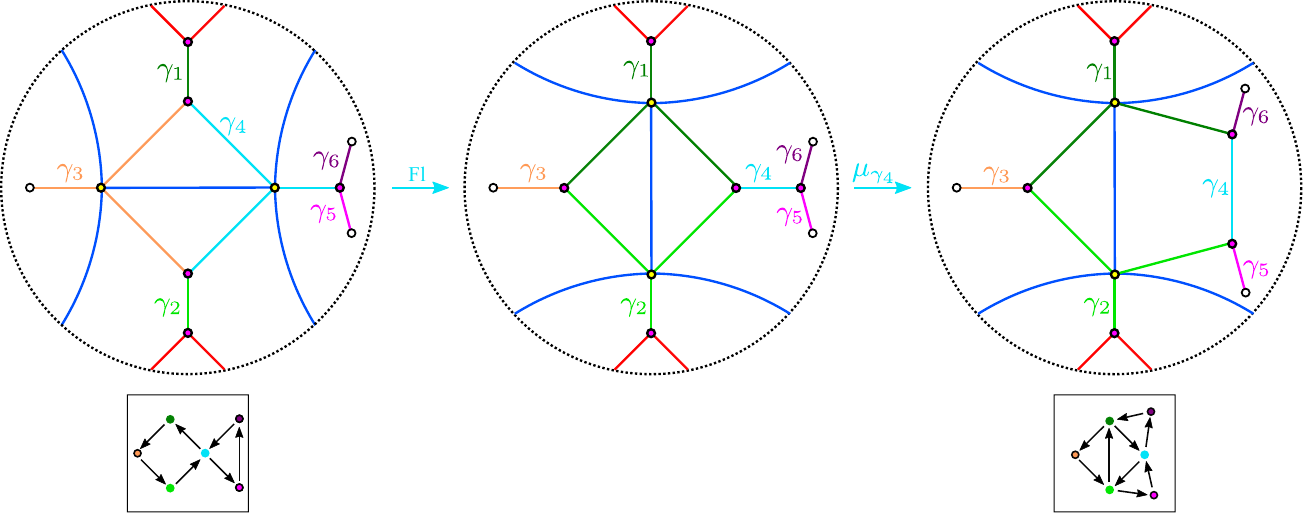}}\caption{Type III.i to Type IV mutations.}\label{fig:3.1to4}
				\end{figure}
			\end{center}
			\item[iii.] (Type III.i to Type IV) We now consider the sequence of 3-graphs in Figure \ref{fig:3.1to4}. Since the initial 3-graph is not sharp at $\gamma_4$, we must first apply a flop before mutating. After applying this flop, $\gamma_4$ is a short {\sf I}-cycle and the 3-graph is sharp at $\gamma_4$. Mutating at $\gamma_4$ then yields a Type IV.i normal form. The short {\sf I}-cycles $\gamma_5$ and $\gamma_6$ are included to indicate where any $A_n$ tail cycles would be sent under this mutation.\\
			
			\begin{center} 
				\begin{figure}[h!]
					{\includegraphics[width=\textwidth]{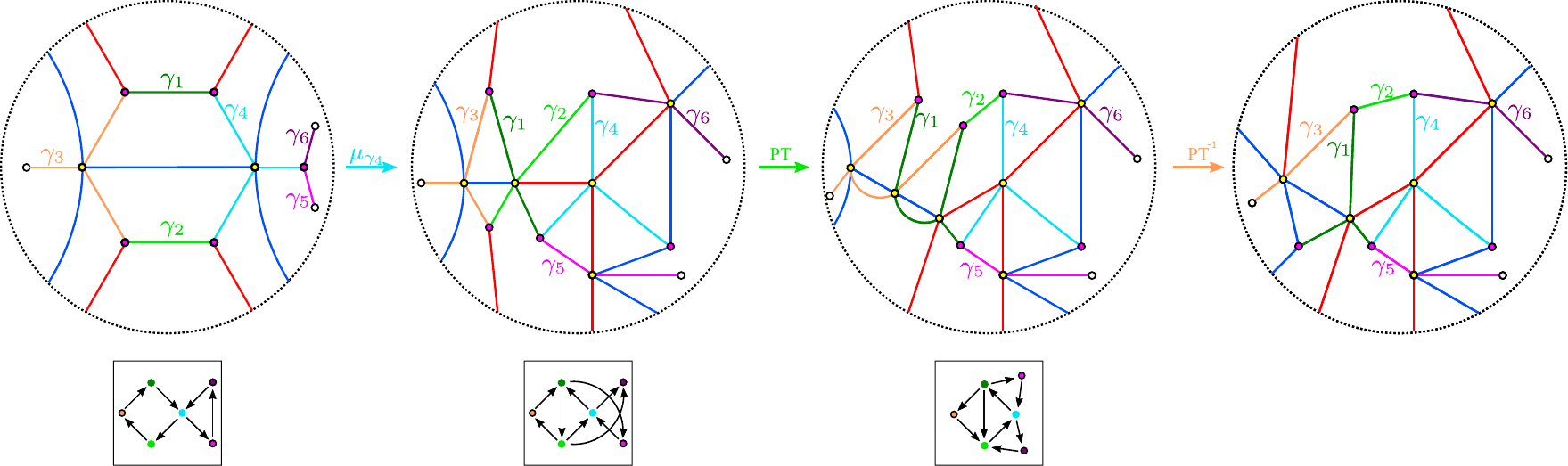}}\caption{Type III.ii to Type IV mutations.}\label{fig:3.2to4}
				\end{figure}
			\end{center}

			\item[iv.] (Type III.ii to Type IV) In Figure \ref{fig:3.2to4},
			mutation at $\gamma_4$ causes $\gamma_1$ and $\gamma_2$ to cross while still intersecting $\gamma_3$ and $\gamma_4$ at either end. We resolve this by first applying a push-through to $\gamma_2$ and then applying a reverse push-through to the trivalent vertex where $\gamma_1$ and $\gamma_3$ intersect a red edge. This results in a sharp 3-graph with $\gamma_1,$ $\gamma_2$, $\gamma_3$, and $\gamma_4$ making up the Type IV normal form. 	We again include $\gamma_5$ and $\gamma_6$ as cycles belonging to a potential $A_n$ tail subgraph in order to show where the $A_n$ tail cycles are sent under this mutation. \\
			
		\end{itemize}

		\textbf{Type IV:} Figure \ref{fig:4Quivers} illustrates all of the relevant Type IV quivers and their mutations. In general, the edges of a Type IV quiver have the form of a single $k-$cycle with the possible existence of 3-cycles or outward-pointing ``spikes" at any of the edges along the $k-$cycle. At the tip of each of these spikes is a possible $A_n$ tail subquiver. We will refer to a vertex at the tip of any of the spikes (e.g., the vertex $v_3$ in Figure \ref{fig:4Quivers}) as a spike vertex and any vertex along the $k-$cycle will be referred to as a $k-$cycle vertex. A homology cycle corresponding to a spike vertex will be referred to as a spike cycle.  Mutating at a spike vertex increases the length of the internal $k-$cycle by one, while mutating at a $k-$cycle vertex decreases the length by 1, so long as $k>3$. Figures \ref{fig:4.1to1}, \ref{fig:4.2to1}, \ref{fig:4.1to3}, and \ref{fig:4.2to3} illustrate the corresponding mutations of 3-graphs for Type IV to Type I and Type IV to Type III when $k=3$. 
		
		\begin{center}
			\begin{figure}[h!]
				{\includegraphics[width=.8\textwidth]{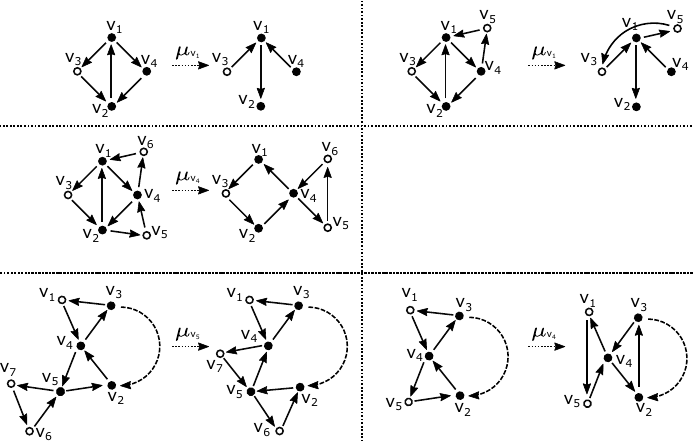}}\caption{From top to bottom, Type IV to Type I, Type IV to Type III, Type IV spike vertex (left) and cycle vertex (right) quiver mutations. The presence or absence of the $A_n$ tail vertices $v_6$ and $v_7$ in the quiver mutation depicted in the first column, third row correspond to the presence or absence of spikes appearing in the resulting quiver.}\label{fig:4Quivers}
			\end{figure}
		\end{center}
		
		\begin{center} 
			\begin{figure}[h!]
				{\includegraphics[width=.9\textwidth]{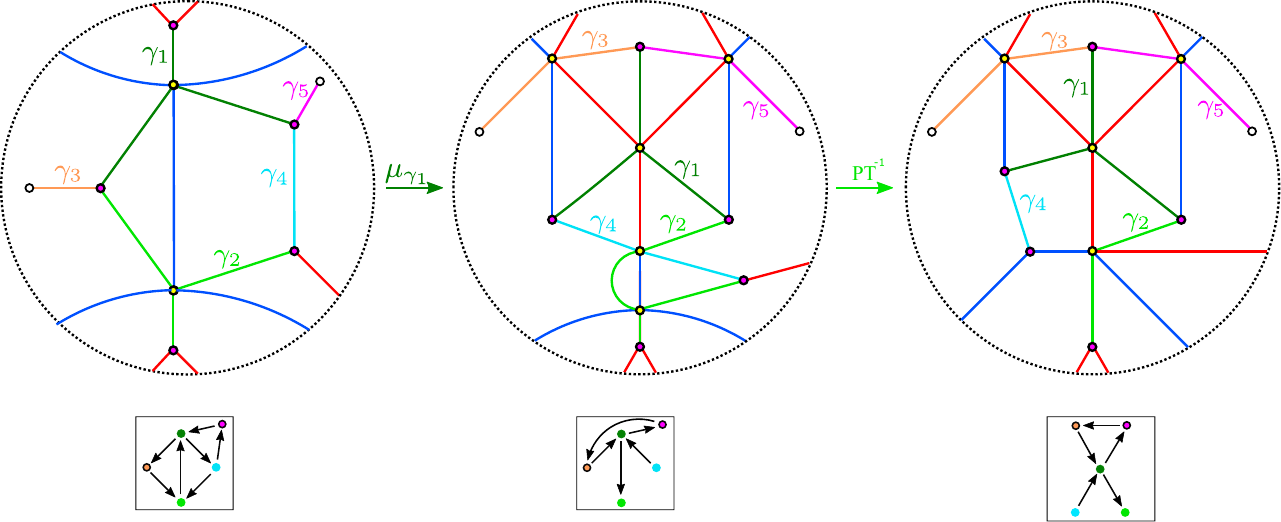}}\caption{Type IV.i to Type I mutations.}\label{fig:4.1to1}
			\end{figure}
		\end{center}

		\begin{itemize}
			\item[i.] (Type IV.i to Type I) We first consider the sequence of 3-graphs in Figure \ref{fig:4.1to1}. Mutation at $\gamma_1$ causes $\gamma_2$ and $\gamma_4$ to cross. Application of a reverse push-through at the trivalent vertex where $\gamma_2$ and $\gamma_4$ intersect a red edge removes this crossing and yields a Type I normal form where $\gamma_1$ is the sole {\sf Y}-cycle.\\
			
			\begin{center} 
				\begin{figure}[h!]
					{\includegraphics[width=.9\textwidth]{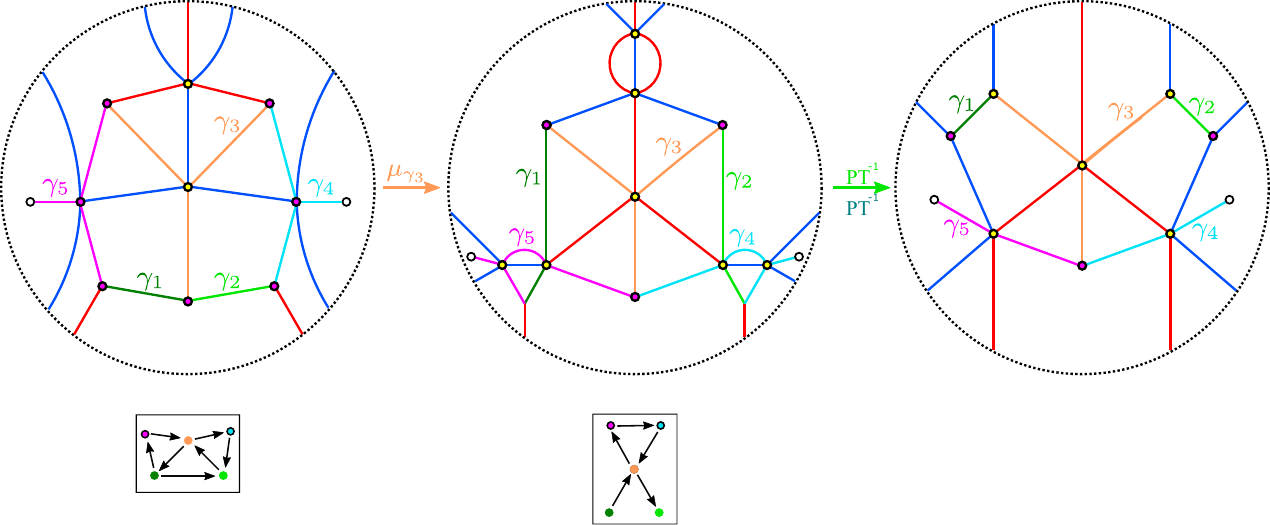}}\caption{Type IV.ii to Type I mutations.}\label{fig:4.2to1}
				\end{figure}
			\end{center}    
			
			\item[ii.] (Type IV.ii to Type I) Mutation at $\gamma_3$ in Figure \ref{fig:4.2to1} yields a 3-graph with geometric intersections between $\gamma_1$ and $\gamma_5$ and between $\gamma_2$ and $\gamma_4$. The application of reverse push-throughs at the trivalent vertex intersections of $\gamma_1$ with $\gamma_5$ and $\gamma_2$ with $\gamma_4$ 
			removes these geometric intersections, resulting in a Type I normal form where $\gamma_1$ is the sole {\sf Y}-cycle. We also apply a candy twist (Legendrian Surface Reidemeister move I) to simplify the intersection at the top of the resulting 3-graph.\\

			\begin{center} 
				\begin{figure}[h!]
					{\includegraphics[width=\textwidth]{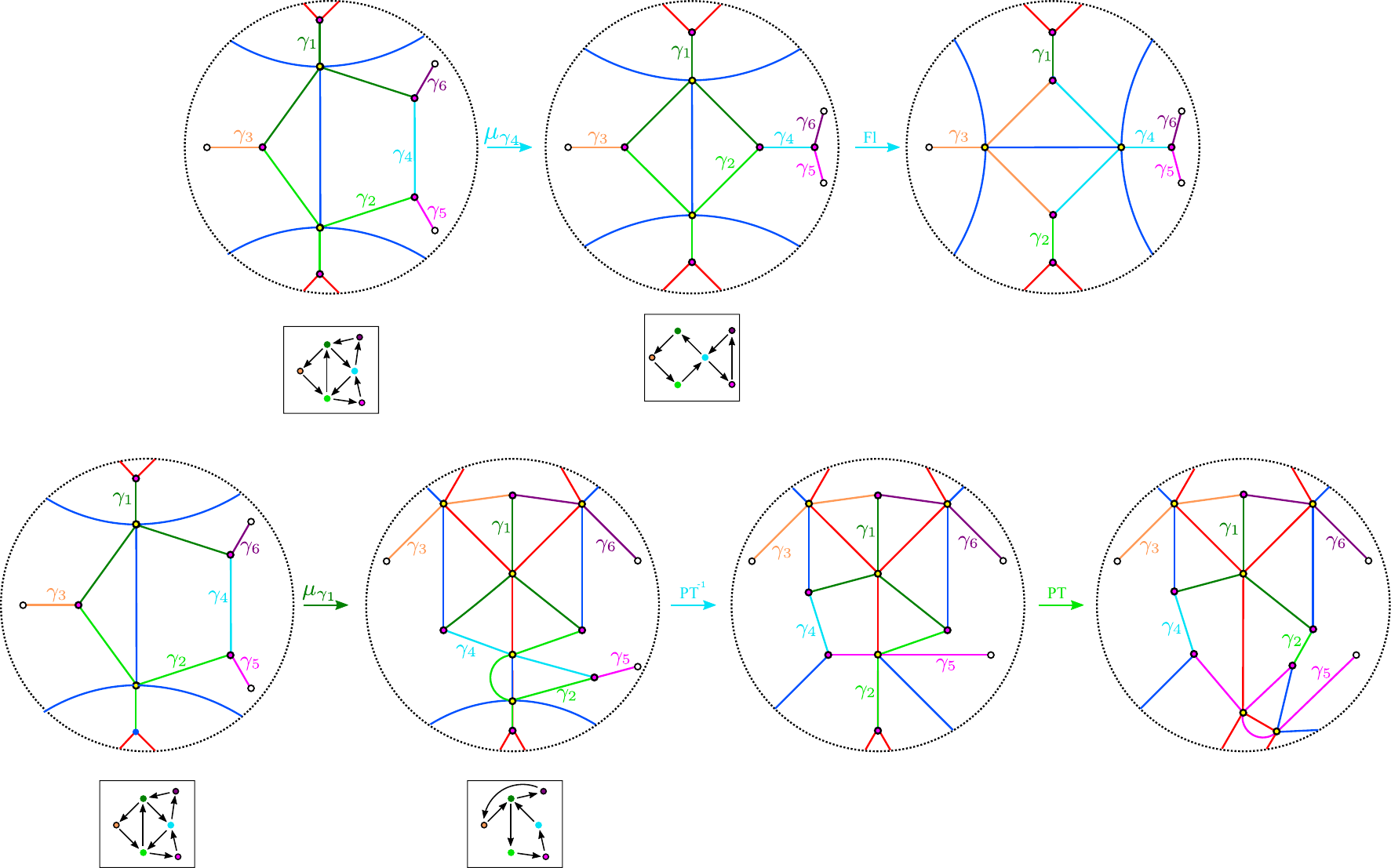}}\caption{Type IV.i to Type III mutations.}\label{fig:4.1to3}
				\end{figure}
			\end{center}

			\item[iii.] (Type IV.i to Type III) We now consider the two sequences of 3-graphs in Figure \ref{fig:4.1to3}. Mutation at any of $\gamma_1, \gamma_2$, $\gamma_3$, or $\gamma_4$ in the Type IV.i normal form yields a Type III normal form. Specifically, mutation at $\gamma_4$ yields a Type III.i normal form that requires no simplification, while mutation at $\gamma_3$ (not pictured) yields a Type III.ii normal form that also requires no simplification. The computation for mutation at $\gamma_1$ is pictured in the sequence on the right and is identical to the computation for mutation at $\gamma_2.$ The first step of the simplification is the same as the Type IV.i to Type I subcase described above. However, we require the application of an additional push-through to remove the geometric intersection between $\gamma_2$ and $\gamma_5.$ This makes $\gamma_5$ into a {\sf Y}-cycle and results in a Type III normal form.\\ 

			\begin{center} 
				\begin{figure}[h!]
					{\includegraphics[width=.9\textwidth]{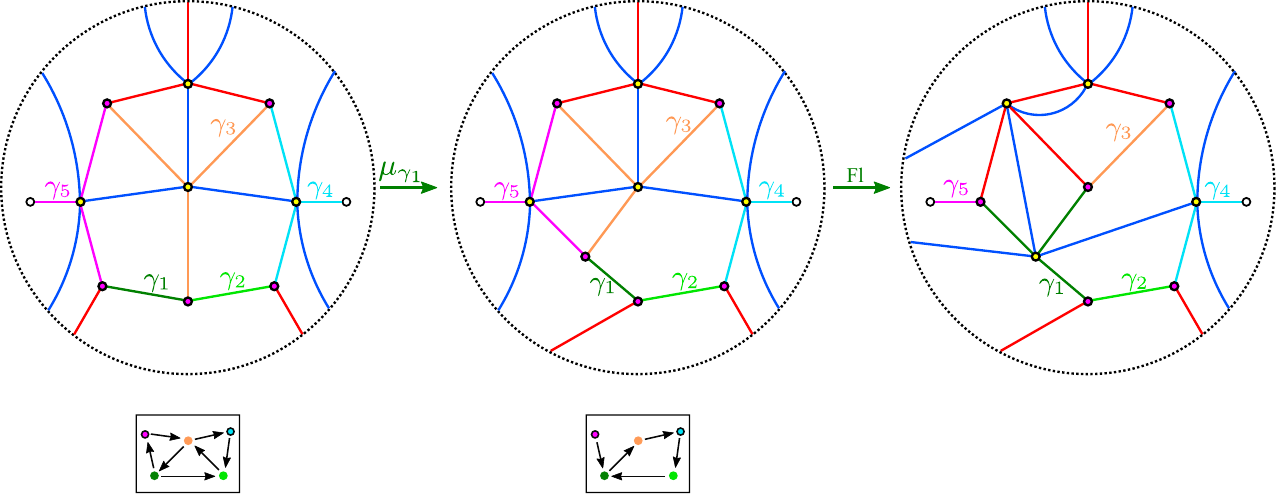}}\caption{Type IV.ii to Type III mutations.}\label{fig:4.2to3}
				\end{figure}
			\end{center}
			
			\item[iv.] (Type IV.ii to Type III)  
			Mutation at $\gamma_1$ in our Type IV.ii normal form, depicted in Figure \ref{fig:4.2to3}, results in a pair of geometric intersections between $\gamma_3$ and $\gamma_5$. Application of a flop removes these geometric intersections and results in a sharp 3-graph with {\sf Y}-cycles $\gamma_1$ and $\gamma_4$, which matches our Type III.ii normal form. Note that the computations for mutations involving a Type IV.ii 3-graph with a single spike cycle are identical.
			
			\bigskip
			The remaining three subcases are all Type IV to Type IV mutations.
			
			\bigskip
			\item[v.] (Type IV.ii to Type IV) Figure \ref{fig:IV.iiSpikes} depicts mutation of a Type IV.ii normal form at a spike cycle.
			Mutating at $\gamma_5$ results in an additional geometric intersection between $\gamma_1$ and $\gamma_3$. We first apply a reverse push-through at the trivalent vertex where $\gamma_1, \gamma_2$ and $\gamma_3$ meet. This introduces an additional geometric intersection between $\gamma_2$ and $\gamma_3$, that we resolve by applying a push-through to $\gamma_3$. Application of a reverse push-through to the trivalent vertex where $\gamma_1$ and $\gamma_5$ intersect a red edge resolves the final geometric intersection between $\gamma_1$ and $\gamma_5$. The {\sf Y}-cycles of the resulting 3-graph correspond to $k-$cycle vertices of the quiver. As shown below, none of the other Type IV to Type IV mutations result in {\sf Y}-cycles corresponding to spike vertices. Therefore, assuming we have simplified after each of our mutations in the manner described above, the only possible way a Type IV.ii 3-graph arises is by mutating from the initial Type I graphs in Figure \ref{fig:1to4.2}. Hence, all other Type IV 3-graphs only have {\sf Y}-cycles corresponding to $k-$cycle vertices in the quiver.	The computations involving a Type IV.ii 3-graph with a single spike cycle are again identical.\\
			
			\begin{center}
				\begin{figure}[h!]
					{\includegraphics[width=\textwidth]{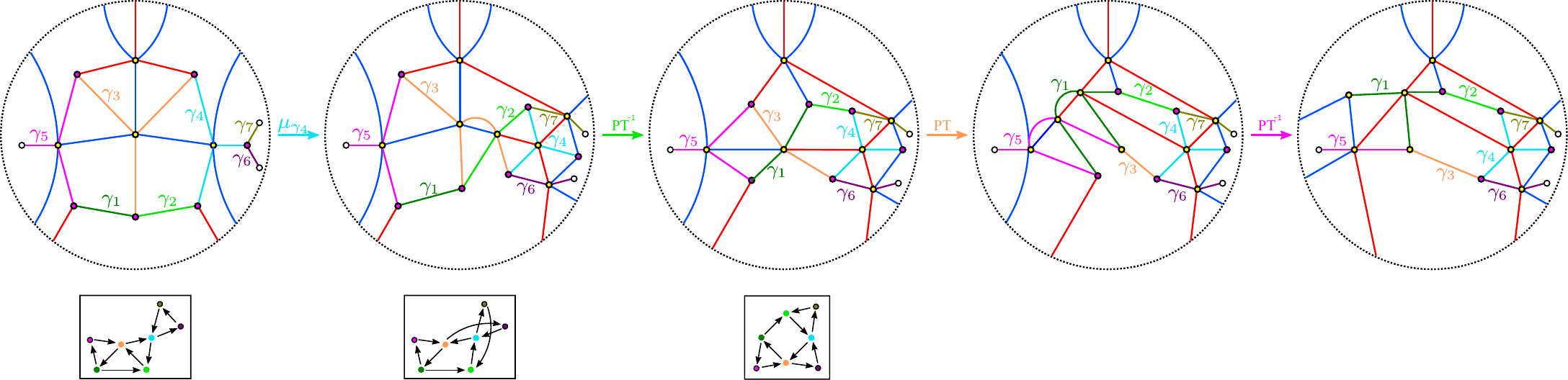}}\caption{Type IV.ii graph mutation at a spike cycle.}\label{fig:IV.iiSpikes}
				\end{figure}
			\end{center}
			
			\item[vi.] (Type IV to Type IV) Figure \ref{fig:4to4k} depicts Type IV to Type IV mutations when the length of the quiver $k-$cycle is greater than 3. When mutating at a homology cycle corresponding to a $k-$cycle vertex of the quiver, we have two possibilities. Figure \ref{fig:4to4k} (top) shows the case where $\gamma_4$ intersects another {\sf Y}-cycle $\gamma_2$, which corresponds to a $k-$cycle vertex in the quiver. Figure \ref{fig:4to4k} (bottom) considers the case where $\gamma_4$ only intersects {\sf I}-cycles. In both of these cases we must apply a reverse push-through to the trivalent vertex where $\gamma_3$ and $\gamma_4$ intersect a red edge in order to simplify the 3-graph. This particular simplification requires that neither of the two edges adjacent to the leftmost edge of $\gamma_4$ carry a cycle before we mutate. A similar computation (not pictured) involving the {\sf Y}-cycle $\gamma_2$ would also require that neither of the two edges adjacent to the bottommost edge of $\gamma_2$ carry a cycle. Crucially, our computations show that Type IV to Type IV mutation preserve this property, i.e., that both of the {\sf Y}-cycles have an edge that is adjacent to a pair of edges which do not carry a cycle. When $k=4,$ the resulting 3-graph resulting from the computations in the top line will have a short {\sf I}-cycle adjacent to $\gamma_2$ and $\gamma_3$, while the 3-graph resulting from the computations in the bottom line will have a short {\sf Y}-cycle adjacent to $\gamma_2$ and $\gamma_3$.\\

			\begin{center}
				\begin{figure}[h!]
					{\includegraphics[width=\textwidth]{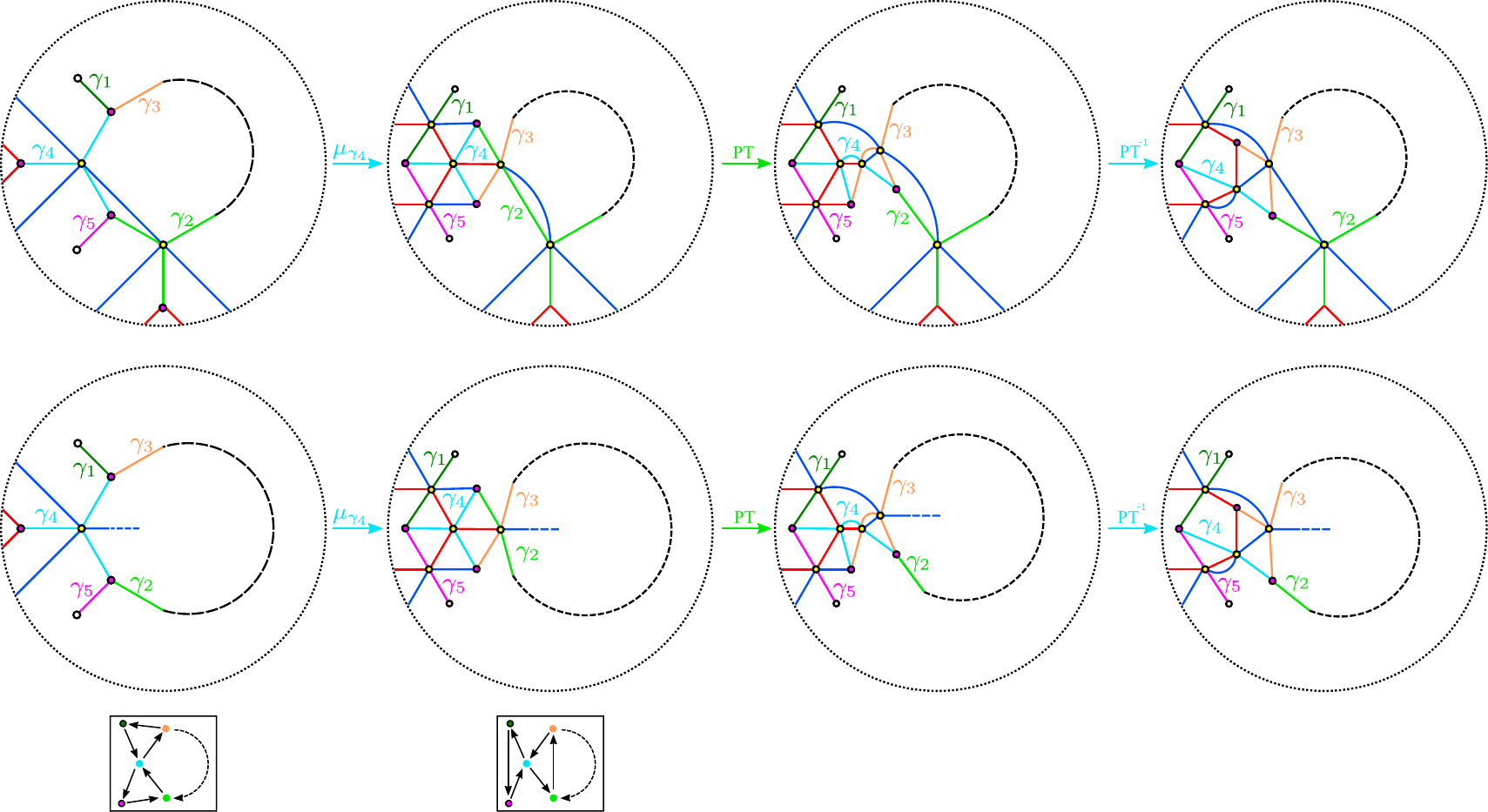}}\caption{Type IV to Type IV mutations at homology cycles corresponding to $k-$cycle vertices in the quiver. Mutating at $\gamma_2,\gamma_3,$ or $\gamma_4$ (corresponding to $k-$cycle vertices in the quiver) in the 3-graphs on the left decreases the length of the $k-$cycle in the quiver by 1.}\label{fig:4to4k}
				\end{figure}
			\end{center}
			
			\begin{center}
				\begin{figure}[h!]
					{\includegraphics[width=.5\textwidth]{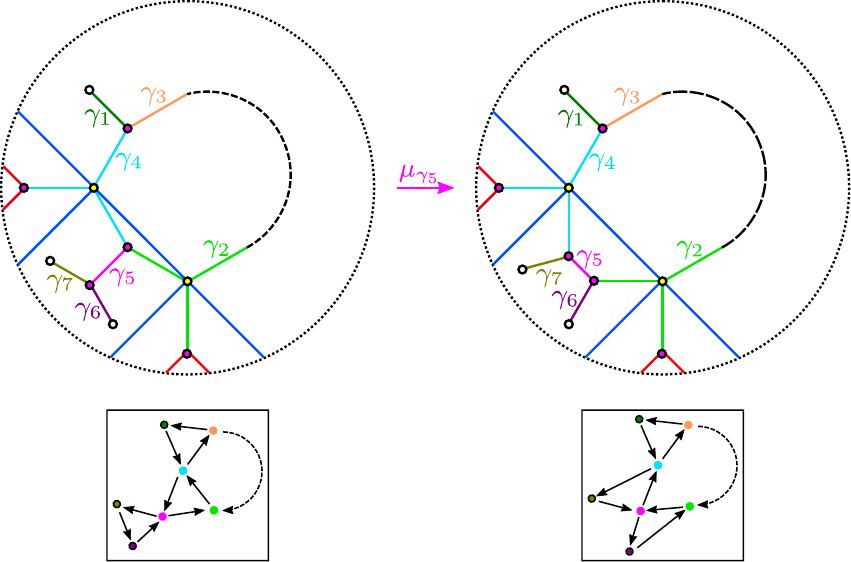}}\caption{Type IV to Type IV mutations at spike cycles. Mutating at the spike cycles $\gamma_1$ or $\gamma_5$ in the 3-graphs on the left increases the length of the $k-$cycle in the intersection quiver by 1.}\label{fig:4to4Spike}
				\end{figure}
			\end{center}
			
			\item[vii.] (Type IV to Type IV) Figure \ref{fig:4to4Spike} depicts mutation at a spike cycle. Since we have already discussed the Type IV.ii spike cycle subcase above, we need only consider the case where each of the spike cycles is a short {\sf I}-cycle. $\gamma_7$ and $\gamma_6$ are included to help indicate where $A_n$ tail cycles are sent under this mutation. The computation for mutating at a spike edge for Type IV.i (i.e. the $k=3$ case) is identical to the $k>3$ case. We have omitted the case where each of the cycles involved in our mutation is an {\sf I}-cycle, but the computation is again a straightforward mutation of a single {\sf I}-cycle that requires no simplification. \\
			

			
		\end{itemize}

		In each of the Type IV to Type IV subcases above, mutating at a {\sf Y}-cycle or an {\sf I}-cycle and applying the simplifications as shown preserves the number of {\sf Y}-cycles in our graph. Therefore, our computations match the normal form we gave in Figure \ref{fig:normalforms} with $k-2$ short {\sf I}-cycles in the normal form 3-graph not belonging to any $A_n$ tail subgraphs.

		This completes our classification of the mutations of normal forms. In each case, we have produced a 3-graph of the correct normal form that is locally sharp and made up of short {\sf Y}-cycles and {\sf I}-cycles. Thus, any sequence of quiver mutations for the intersection quiver $Q(\Gamma_0(D_n),\{\gamma_i^{(0)}\})$ of our initial $\Gamma_0(D_n)$ is weave realizable. Hence, given any sequence of quiver mutations, we can apply a sequence of Legendrian mutations to our original 3-graph to arrive at a 3-graph with intersection quiver given by applying that sequence of quiver mutations to $Q(\Gamma_0(D_n),\{\gamma_i^{(0)}\})$, as desired.

	\end{proof}
	
	Having proven weave realizability for $\Gamma_0(D_n)$, we conclude with a proof of Corollary \ref{cor}.
	
	\subsection{Proof of Corollary \ref{cor}}
	We take our initial cluster seed in $\mathcal{C}(\Gamma)$ to be the cluster seed associated to $\Gamma_0(D_n)$. The cluster variables in this initial seed exactly correspond to the microlocal monodromies along each of the homology cycles of the initial basis $\{\gamma_i^{(0)}\}$. The intersection quiver $Q(\Gamma_0(D_n), \{\gamma_i^{0}\})$ is the $D_n$ Dynkin diagram and thus the cluster seed is $D_n$-type. By definition, any other cluster seed in the $D_n$-type cluster algebra is obtained by a sequence of quiver mutations starting with the quiver $Q(\Gamma_0(D_n), \{\gamma_i^{0}\})$ and its associated cluster variables. Theorem \ref{main} implies that any quiver mutation of $Q(\Gamma_0(D_n), \{\gamma_i^{0}\})$ can be realized by a Legendrian mutation in $\Lambda(\Gamma_0(D_n)),$  so we have proven the first part of the corollary. The remaining part of the corollary follows from the fact that the $D_n$-type cluster algebra is known to be of finite mutation type with $(3n-2)C_{n-1}$ distinct cluster seeds. \hfill $\Box$


	\subsection{Further Study} While a classification of $E$-type quivers is not yet known, it seems likely that the techniques in this manuscript could be used to show weave realizability for Lagrangian fillings arising from $\lambda(E_6), \lambda(E_7),$ and $\la(E_8).$ Identifying normal forms for the expected weave fillings \cite[Conjecture 5.1]{CasalsLagSkel} could even aid in such a classification of $E$-type quivers. More generally, it is possible that the methods used here may be adapted to show weave realizability for any positive braid.


	
	\bibliographystyle{alpha}
	\bibliography{DnType}

\end{document}